\documentclass{article}
\usepackage{amsfonts, amsmath, amstext, amssymb}

\usepackage[dvips]{graphicx}
\usepackage{latexsym}

\newtheorem{thm}{Theorem}[section]
\newtheorem{prop}[thm]{Proposition}

\newtheorem{lemma}[thm]{Lemma}

\newtheorem{cor}[thm]{Corollary}

\newtheorem{claim}[thm]{Claim}

\newtheorem{defn}[thm]{Definition}

\newenvironment{pf}{\begin{trivlist}\item[\hskip\labelsep
{\it Proof.}]}{\end{trivlist}}

\newcommand{\Cof}{\text{Cof}}
\newcommand{\Fin}{\text{Fin}}

\newcommand{\set}[2]{\ensuremath{ \{ #1 : #2 \} }}

\newcommand{\Z}{\mathbb{Z}}
\newcommand{\Q}{\mathbb{Q}}

\newcommand{\E}{\mathcal{E}}

\renewcommand{\S}{\mathcal{S}}

\def\converges{\!\downarrow}
\newcommand{\at}{\char'100}

\newcommand{\la}{\langle}
\newcommand{\ra}{\rangle}

\newcommand{\qed}{\hbox to 0pt{}\nobreak\hfill\rule{2mm}{2mm}}

\newcommand{\xvec}{\vec{x}}
\newcommand{\bfz}{\boldsymbol{0}}

\newcommand{\bfd}{\boldsymbol{d}}

\newcommand{\dom}[1]{\text{dom}(#1)}

\newcommand{\EIso}[1]{E^{\text{#1}}_{\cong}}
\newcommand{\FIso}{F^{\text{alg}}_{\cong}}
\newcommand{\FAC}{F_{\cong}^{AC}}
\newcommand{\sqrtpn}{\sqrt{p_n}}

\def\Qbar{\overline{\mathbb{Q}}}

\def\s01{\ensuremath{\Sigma^0_1}}
\def\d02{\ensuremath{\Delta^0_2}}
\def\phi{\varphi}
\def\ahat{\ensuremath{\hat{a}}}
\def\bhat{\ensuremath{\hat{b}}}
\def\chat{\ensuremath{\hat{c}}}
\def\dhat{\ensuremath{\hat{d}}}
\def\xhat{\ensuremath{\hat{x}}}

\def\xvec{\ensuremath{\vec{x}}}
\def\yvec{\ensuremath{\vec{y}}}
\def\Ahat{\ensuremath{\widehat{A}}}
\def\mhat{\ensuremath{\hat{m}}}
\def\What{\ensuremath{\widehat{W}}}
\def\res{\!\!\upharpoonright\!}

\def\Vec{\textbf{Vec}}

\newcommand{\Ece}[1]{\ensuremath{E^{ce}_{#1}}}
\newcommand{\Ecemax}{\ensuremath{E^{ce}_{\text{max}}}}
\newcommand{\EDmax}{\ensuremath{E^{D}_{\text{max}}}}
\newcommand{\Ecemin}{\ensuremath{E^{ce}_{\text{min}}}}
\newcommand{\Eperm}{\ensuremath{E_{\text{perm}}}}
\newcommand{\Eceperm}{\ensuremath{E^{ce}_{\text{perm}}}}
\newcommand{\Eceset}{\ensuremath{E^{ce}_{\text{set}}}}
\newcommand{\Ececard}{\ensuremath{E^{ce}_{\text{card}}}}
\newcommand{\EDcard}{\ensuremath{E^{D}_{\text{card}}}}
\newcommand{\Ecardprime}{\ensuremath{E^{\emptyset'}_{\text{card}}}}
\newcommand{\comment}[1]{}
\newcommand{\leqset}{\leq_{set}}
\newcommand{\geqset}{\geq_{set}}

\newcommand{\Zce}{\ensuremath{Z^{ce}_{0}}}
\newcommand{\rs}[1]{\upharpoonright{#1}}

\newcommand{\ER}{equivalence relation}
\newcommand{\ERs}{equivalence relations}

\begin{document}

\title{Finitary reducibility on equivalence relations}
\author{Russell Miller~\&~Keng Meng Ng
%\thanks{The first author was
%partially supported by grants number 67182-00-36,
%68470-00 37, and 80209-04-12 from
%The City University of New York PSC-CUNY Research Award Program.
%The collaboration of these authors was facilitated by the program
%\emph{Sets and Computations} at the Institute for Mathematical Sciences
%of the National University of Singapore in April 2015.}
}

\maketitle

\begin{abstract}
We introduce the notion of finitary computable reducibility
on equivalence relations on the domain $\omega$.
This is a weakening of the usual notion of computable reducibility,
and we show it to be distinct in several ways.  In particular,
whereas no equivalence relation can be $\Pi^0_{n+2}$-complete
under computable reducibility, we show that,
for every $n$, there does exist a natural equivalence relation
which is $\Pi^0_{n+2}$-complete under finitary reducibility.
We also show that our hierarchy of finitary reducibilities does not
collapse, and illustrate how it sharpens certain known results.
Along the way, we present several
new results which use computable reducibility to establish
the complexity of various naturally defined equivalence
relations in the arithmetical hierarchy.  
%MYHILL  We also refute a possible generalization of Myhill's Theorem.
\end{abstract}

\section{Introduction}
\label{sec:intro}

Computable reducibility provides a natural way of measuring
and comparing the complexity of equivalence relations
on the natural numbers.  Like most notions of reducibility
on sets of natural numbers, it relies on the concept of
Turing computability to rank objects according to their complexity,
even when those objects themselves may be far from computable.
It has found particular usefulness in computable model theory,
as a measurement of the classical property of being isomorphic:
if one can computably reduce the isomorphism problem
for computable models of a theory $T_0$ to the isomorphism
problem for computable models of another theory $T_1$,
then it is reasonable to say that isomorphism on models of $T_0$
is no more difficult than on models of $T_1$.  The related notion
of Borel reducibility was famously applied this way by Friedman
and Stanley in \cite{FS89}, to study the isomorphism problem
on all countable models of a theory.
Yet computable reducibility has also become the subject
of study in pure computability theory, as a way of ranking
various well-known equivalence relations arising there.

Recently, as part of our study of this topic,
we came to consider certain reducibilities weaker than
computable reducibility.  This article introduces
these new, finitary notions of reducibility on equivalence
relations and explains some of their uses.  We believe that
researchers familiar with computable reducibility will find
finitary reducibility to be a natural and appropriate measure
of complexity, not to supplant computable reducibility but to
enhance it and provide a finer analysis of situations in which
computable reducibility fails to hold.

Computable reducibility is readily defined.  It has gone by many
different names in the literature, having been called $m$-reducibility
in \cite{ceers,BS83,GG01} and FF-reducibility in \cite{FF09,FF11,FFT10},
in addition to a version on first-order
theories which was called Turing-computable reducibility
(see \cite{CCKM04,CK06}).
\begin{defn}
\label{defn:compreducibility}
Let $E$ and $F$ be equivalence relations on $\omega$.
A \emph{reduction} from $E$ to $F$ is a function
$g:\omega\to\omega$ such that
\begin{equation}
\label{eq:reduction}
 \forall x,y\in\omega~~[x~E~y%\qquad
 ~~\iff%\qquad 
 ~~g(x)~F~g(y)].
 \end{equation}
We say that $E$ is \emph{computably reducible} to $F$,
written $E\leq_c F$, if there exists a reduction from $E$
to $F$ which is Turing-computable.  More generally,
for any Turing degree $\bfd$, $E$ is
$\bfd$-\emph{computably reducible} to $F$
if there exists a reduction from $E$
to $F$ which is $\bfd$-computable.
\end{defn}
There is a close analogy between this definition and that of
\emph{Borel reducibility}:  in the latter, one considers equivalence
relations $E$ and $F$ on the set $2^\omega$ of real numbers,
and requires that the reduction $g$
be a Borel function on $2^\omega$.
In another variant, one requires $g$ to be
a continuous function on reals (i.e., given by a Turing functional
$\Phi^Z$ with an arbitrary real oracle $Z$), thus defining
\emph{continuous reducibility} on equivalence relations on $2^\omega$.
\comment{
%%%%
For relations on $\omega$, it is useful to generalize Definition
\ref{defn:compreducibility} slightly to the following.
\begin{defn}
\label{defn:dcompreducibility}
Let $E$ and $F$ be equivalence relations on $\omega$,
and $\bfd$ any Turing degree.
We say that $E$ is \emph{$\bfd$-computably reducible} to $F$,
written $E\leq_{\bfd} F$, if there exists a total $\bfd$-computable function $g$
such that
$$ \forall x,y\in\omega~[x~E~y\qquad\iff\qquad g(x)~F~g(y)].$$
\end{defn}
So $\bfz$-computable reducibility is just the notion defined above.
%%%%
}

So a reduction from $E$ to $F$ maps every element in the field of the relation $E$
to some element in the field of $F$, respecting these equivalence relations.
Our new notions begin with \emph{binary computable reducibility}.
In some situations, while it is not possible to give a computable reduction
from $E$ to $F$, there does exist a computable function which takes each pair
$\la x_0,x_1\ra $ of elements from the field of $E$ and outputs a pair
of elements $\la y_0,y_1\ra$ from that of $F$ such that $y_0Fy_1$ if and only if
$x_0Ex_1$.  (The reader may notice that this is simply
an $m$-reduction from the set $E$ to the set $F$.)
Likewise, an \emph{$n$-ary computable reduction} accepts
$n$-tuples $\xvec$ from the field of $E$ and outputs $n$-tuples
$\yvec$ from $F$ with $(x_iEx_j\iff y_iFy_j)$ for all $i<j<n$,
and a \emph{finitary computable reduction} does the same for all
finite tuples.  Intuitively, a computable reduction (as in Definition
\ref{defn:compreducibility}) does the same for all elements from the
field of $E$ simultaneously.

A computable reduction clearly gives us a
computable finitary reduction, and hence a computable $n$-reduction
for every $n$.  Oftentimes, when one builds a computable reduction, 
one attempts the opposite procedure:  the first step is to build
a binary reduction, and if this is successful, one then treats
the binary reduction as a basic module and attempts to combine countably
many basic modules into a single effective construction.  Our initial
encounter with finitary reducibility arose when we found a basic module
of this sort, but realized that it was only possible to combine finitely
many such modules together effectively.

At first we did not expect much from this new notion,
but we found it to be of increasing interest as we continued
to examine it.  For example, we found that the standard $\Pi^0_{n+2}$
equivalence relation defined by equality of the sets $W_i^{\emptyset^{(n)}}$
and $W_j^{\emptyset^{(n)}}$ is complete among $\Pi^0_{n+2}$
equivalence relations under finitary reducibility.  This is of particular
interest because, for precisely these classes, no equivalence relation
can be complete under computable reducibility (as shown recently
in \cite{IMNNS13}).  Extending our study to certain relations from
computable model theory, we found that the isomorphism problem
$\FAC$ for computable algebraically closed fields of characteristic $0$,
while $\Pi^0_3$-complete as a set,
fails to be complete under finitary reducibility:  it is complete
for $3$-ary reducibility, but not for the $4$-ary version.  This confirms
one's intuition that isomorphism on algebraically closed fields, despite
being $\Pi^0_3$-complete as a set, is not
an especially difficult problem, requiring only knowledge of the
transcendence degree of the field.  In contrast, the isomorphism problem
$\FIso$ for algebraic fields of characteristic $0$, while only $\Pi^0_2$,
does turn out to be complete at that level under finitary reducibility.
%(which is as much as one could ask, there being no $\Pi^0_2$-complete
%equivalence relation under computable reducibility).

This paper proceeds much as our investigations proceeded.
In Section \ref{sec:others} we present the equivalence relations
on $\omega$ which we set out to study.  We derive a number of results about them,
and by the time we reach Proposition \ref{prop:binaryEsetE3},
it should seem clear to the reader how the notion of finitary
reducibility arose for us, and why it seems natural in this context.
The exact definitions of $n$-ary and finitary reducibility appear
as Definition \ref{defn:narily}.  In Sections \ref{sec:introfinitary} and
\ref{sec:finitary}, we study finitary reducibility in its own right.
We produce the natural $\Pi^0_{n+2}$ \ERs\ described above,
defined by equality among $\Sigma^0_n$ sets, which are complete
under finitary reducibility among all $\Pi^0_{n+2}$ \ERs.
%, a result of particular
%interest since it is known that, precisely when $m\geq 2$,
%no \ER\ can be $\Pi^0_m$-complete under computable reducibility.
Subsequently we show that the hierarchy of $n$-ary reducibilities
does not collapse, and that several standard equivalence relations
on $\omega$ witness this non-collapse for certain $n$.

\comment{MYHILL

Finally, in Section \ref{sec:Myhill}, we establish some further results
on computable reducibility, including a proof that Myhill's Theorem
does not apply to the relation of computable reducibility,
even in a very simple context.
}

\section{Background in Computable Reducibility}
\label{sec:others}

The purpose of this section is twofold.  First, for the reader who is not
already familiar with the framework and standard methods used
in its study, it introduces some examples of results in computable
reducibility, with proofs.    The examples, however, are not intended
as a broad outline of the subject; they are confined to one very specific
subclass of equivalence relations (those which, as sets, are $\Pi^0_4$),
rather than offering a survey of important results in the field.
In fact the results we prove here are new, to our knowledge.
They use computable reducibility to establish
the complexity of various naturally defined equivalence
relations in the arithmetical hierarchy.  In doing so,
we continue the program of work already set in motion in
\cite{E77,BS83,GG01,CHM12,ceers,IMNNS13} and augment their results.
However, the second and more important purpose of these results is to
help explain how we came to develop the notion of finitary
reducibility and why we find that notion to be both natural
and useful.  By the end of the section, the reader will have an
informal understanding of finitary reducibility, which is then
formally defined and explored in the ensuing two sections.

The following definition introduces several natural
equivalence relations which we will consider in this section.
Here, for a set $A\subseteq\omega$, we write $A^{[n]}=\set{x}{\la x,n\ra\in A}$
for the $n$-th column of $A$ when $\omega$ is viewed as the two-dimensional
array $\omega^2$ under the standard computable pairing function $\la\cdot,\cdot\ra$
from $\omega^2$ onto $\omega$.

\begin{defn}
\label{defn:ERs}
First we define several equivalence relations on $2^\omega$.
\begin{itemize}

\item $E_{perm}= \{\langle A,B\rangle \mid
(\exists $ a permutation $p:\omega\to\omega) (\forall n) A^{[n]}=B^{[p(n)]}\}$.

\item $E_{\Cof} =\{\langle A,B\rangle \mid $ For every $n$, $A^{[n]}$ is
cofinite iff $B^{[n]}$ is cofinite$\}$.

\item $E_{\Fin} =\{\langle A,B\rangle \mid $ For every $n$, $A^{[n]}$ is
finite iff $B^{[n]}$ is finite$\}$.

%\item $E_{=}^n=\{(i,j)\mid W_i^{\emptyset^{(n)}}=W_j^{\emptyset^{(n)}}\}$, for each $n\in\omega$.

%\item $E_{max}^n=\{(i,j)\mid \max W_i^{\emptyset^{(n)}}=\max W_j^{\emptyset^{(n)}}\}$, for each $n\in\omega$.
\end{itemize}
Each of these relations induces an equivalence relation on $\omega$,
by restricting to the c.e.\ subsets of $\omega$ and then allowing the index
$e$ to represent the set $W_e$, under the standard indexing of c.e.\ sets.
The superscript ``ce'' denotes this, so that, for instance,
$$  \Eceperm =\{\langle i,j\rangle \mid (\exists\text{~a permutation~}
p:\omega\to\omega) (\forall n) W_i^{[n]}=W_j^{[p(n)]}\}.$$
Similarly we define $\Ece{\Cof}$ and $\Ece{\Fin}$, and also the
following two equivalence relations on $\omega$ (where the superscripts
denote oracle sets, so that $W_i^D=\dom{\Phi_i^D}$):
\begin{itemize}
\item $E_{=}^n=\{(i,j)\mid W_i^{\emptyset^{(n)}}=W_j^{\emptyset^{(n)}}\}$, for each $n\in\omega$.

\item $E_{max}^n=\{(i,j)\mid \max W_i^{\emptyset^{(n)}}=\max W_j^{\emptyset^{(n)}}\}$, for each $n\in\omega$.
\end{itemize}
In $E_{max}^n$, for any two infinite sets $W_i^{\emptyset^{(n)}}$ and $W_j^{\emptyset^{(n)}}$,
this defines $\la i,j\ra\in E_{max}^n$, since we consider both sets to have
the same maximum $+\infty$.

\comment{
Finally, from computable model theory, we will consider the isomorphism relation on
rational vector spaces,
$$\Vec =\{\langle W,V\rangle \mid W=V\text{~or~}W,V\text{~define
$\Q$-vector spaces such that~}W\cong V\},$$
and its restriction from $2^\omega$ to the world of computable rational vector spaces:
$$\Vec =\{\langle i,j\rangle \mid i=j\text{~or~}W_i,W_j\text{~define
$\Q$-vector spaces such that~}W_i \cong W_j\},$$
(A set $A$, viewed as a subset of $\omega^3$, \emph{defines a rational vector space}
if $\omega$ becomes such a space when $A$ is viewed as an addition relation on $\omega$.)
}
\end{defn}

\subsection{$\Pi^0_4$ equivalence relations}
Here we will clarify the relationship between several equivalence relations
occurring naturally at the $\Pi^0_4$ level.
Recall the equivalence relations $E_3$, $E_{set}$, and $Z_0$
defined in the Borel theory.  Again the analogues of these for c.e.\ sets are
relations on the natural numbers, defined using the symmetric difference $\triangle$:
\begin{align*}
& i~\Ece3~j \qquad\iff\qquad \forall n~[|(W_i)^{[n]} \triangle (W_j)^{[n]} |<\infty]\\
& i~\Eceset~j \qquad\iff\qquad \{(W_i)^{[n]} \mid n\in\omega\}=\{(W_j)^{[n]} \mid n\in\omega\}\\
& i~\Zce~j \qquad\iff\qquad \lim_n \frac{|(W_i\triangle W_j)\rs{n}|}{n}=0
\end{align*}
The aim of this section is to show that the situation in the following picture holds
for computable reducibility.

\begin{center}
\vspace*{0.5cm}
\begin{minipage}{8cm}
\begin{picture}(200,50)
\linethickness{0.6pt}
\put(100,20){\line(0,1){13}}
\put(28,40){$\Eceset\equiv_c\Eceperm\equiv_c\Ece{\Cof}\equiv_c E_=^2$}
\put(76,9){$\Ece{3}\equiv_c \Zce$}
\end{picture}
\end{minipage}
\vspace*{0.5cm}
\end{center}
Hence all these classes fall into two distinct computable-reducibility degrees,
one strictly below the other.  Even though no $\Pi^0_4$ class
is complete under $\leq_c$, we will show that each of these classes
is complete under a more general reduction.

The three classes $\Ece{3},\Ece{\text{set}}$ and $\Zce$ are easily seen to be $\Pi^0_4$.
This is not as obvious for \Eceperm.
\begin{lemma}
\label{lemma:Pi4}
The relation \Eceperm is $\Pi^0_4$, being defined on pairs $\la e,j\ra$ by:
%\begin{align*}
%&\la e,j\ra\in\Eceperm \iff\\&
$$\forall k\forall n_0<\cdots<n_k~
\exists{~distinct~}m_0,\ldots,m_k~\forall i\leq k~(W_e^{[n_i]}
=W_j^{[m_i]}),$$
%\end{align*}
in conjunction with the symmetric statement with $W_j$ and $W_e$
interchanged.
\end{lemma}
\begin{pf}
Since ``$W_e^{[n_i]}=W_j^{[m_i]}$'' is $\Pi^0_2$, the given
statement is $\Pi^0_4$, as is the interchanged version.
The statements clearly hold for all $\la e,j\ra\in\Eceperm$.
Conversely, if the statements hold, then each c.e.\ set which
occurs at least $k$ times as a column in $W_e$
must also occur at least $k$ times as a column in $W_j$,
and vice versa.  It follows that every c.e.\ set occurs equally
many times as a column in each, allowing an easy definition
of the permutation $p$ to show $\la e,j\ra\in\Eceperm$.
\qed\end{pf}

\begin{thm}
\label{thm:EsetEperm}
\Eceperm\ and \Eceset\ are computably bireducible. (We write $\Eceperm\equiv_c \Eceset$ to denote this.)
\end{thm}
\begin{pf}
For the easier direction \Eceset $\leq_c$ \Eceperm, given a c.e. set $A$,
define uniformly the c.e.\ set $\widehat{A}$ by setting (for each $e,i,x$)
$x\in \widehat{A}^{[\langle e,i\rangle]}$ iff $x\in A^{[e]}$. That is, we repeat each column of $A$
infinitely many times in $\widehat{A}$. Then $A~E_{\text{set}}~ B$ iff
$\widehat{A}~E_{\text{perm}}~\widehat{B}$.  (Since the definition is uniform,
there is a computable function $g$ which maps each $i$ with $W_i=A$
to $g(i)$ with $W_{g(i)}=\widehat{A}$.  This $g$ is the computable reduction
required by the theorem, with $i~\Eceset~j$ iff $g(i)~\Eceperm~g(j)$ for all $i,j$.)

We now turn to $\Eceperm\leq_c\Eceset$. Fix a c.e.\ set $A$. We describe
a uniform procedure to build $\widehat{A}$ from $A$. We must do this in a way where for any pair of c.e. sets $W,V$, $W\Eceperm V$ iff $\widehat{W}\Eceset\widehat{V}$. The computable function $q$ that gives $W_{q(i)}=\widehat{W}_i$ will then be a witness for the reduction $\Eceperm\leq_c\Eceset$.

For each $x$ let $F(x)$
be the number of columns $y\leq x$ such that $A^{[x]}=A^{[y]}$. There is a
natural computable guessing function $F_s(x)$ such that for every $s$,
$F_s(x)\leq x$ and $F(x)=\lim\sup_s F_s(x)$.

Associated with $x$ are the c.e.\ sets $C[x,n]$ for each $n>0$ and $D[x,i,j]$ for each $i>0,j\in\omega$, defined as follows. $D[x,i,j]$ is the set $D$ such that
\[         D^{[k]} =
         \begin{cases}
            A^{[x]}, & \text{if $k =0$,}\\
            \{0,1,\cdots,j-1\}, & \text{if $k=i$,}\\
            \emptyset, & \text{otherwise.}\\
         \end{cases}
\]
and $C[x,n]$ is the set $C$ such that
\[         C^{[k]} =
         \begin{cases}
            A^{[x]}, & \text{if $k =0$,}\\
            %\{0,1,\cdots,\max\{s:F_s(x)\geq n\}\}, & \text{if $k=n$, and $\forall^\infty s (F_s(x)< n)$,}\\
            %\omega, & \text{if $k=n$, and $\exists^\infty s (F_s(x)\geq n)$,}\\
            \left\{t:(\exists s\geq t)(F_s(x)\geq n)\right\}, & \text{if $k=n$},\\
						\emptyset, & \text{otherwise.}\\
         \end{cases}
\]
Now let $\widehat{A}$ be obtained by copying all the sets $C[x,n]$ and $D[x,i,j]$
into the columns. That is, let $\widehat{A}^{[2\langle x,n\rangle]}=C[x,n]$ and
$\widehat{A}^{[2\langle x,i,j\rangle+1]}=D[x,i,j]$. Now suppose that $A~\Eperm B$.
We verify that $\widehat{A}~E_{\text{set}} \widehat{B}$, writing $C[A,x,n]$, $C[B,x,n]$, $D[A,x,i,j]$,
and $D[B,x,i,j]$ to distinguish between the columns of $\widehat{A}$ and $\widehat{B}$.

Fix $x$ and consider $D[A,x,i,j]$. Since there is some $y$ such that $A^{[x]}=B^{[y]}$
it follows that $D[A,x,i,j]=D[B,y,i,j]$ for every $i,j$. Now we may pick $y$ such that
$F(A,x)=F(B,y)$. It then follows that $C[A,x,n]=C[B,y,n]$ for every $n\leq F(A,x)$,
and for $n>F(A,x)$ we have $C[A,x,n]=D[B,y,n,j]$ for some appropriate $j$.
Hence every column of $\widehat{A}$ appears as a column of $\widehat{B}$.
A symmetric argument works to show that every column of $\widehat{B}$ is a column
of $\widehat{A}$.

Now suppose that $\widehat{A}~E_{\text{set}}~\widehat{B}$. We argue that $A~E_{\text{perm}}~B$.
Fix $x$ and $n$ such that there are exactly $n$ many different numbers $z\leq x$
with $A^{[z]}=A^{[x]}$. We claim that there is some $y$ such that $A^{[x]}=B^{[y]}$ and there
are at least $n$ many $z\leq y$ such that $B^{[z]}=B^{[y]}$.

The column $C[A,x,n]$ of $\widehat{A}$ is the set $C$ such that $C^{[0]}=A^{[x]}$ and
$C^{[n]}=\omega$. Now $C[A,x,n]$ cannot equal $D[B,y,i,j]$ for any $y,i,j$
since $D$-sets have every column finite except possibly for the $0^{\text{th}}$ column.
So $C[A,x,n]=C[B,y,n]$ for some $y$. It follows that $A^{[x]}=\left(C[B,y,n]\right)^{[0]}=B^{[y]}$,
and we must have $\lim\sup_s F_s(B,y)\geq n$. So each $A^{[x]}$ corresponds
to a column $B^{[{y'}]}$ of $B$ with $F(B,y')= F(A,x)$. Again a symmetric argument
follows to show that each $B^{[y]}$ corresponds to a column $A^{[x]}$ of $A$
with $F(A,x)= F(B,y)$. Hence $A$ and $B$ agree up to a permutation of columns.
\qed
\end{pf}

\begin{thm}
\label{thm:EsetEcof}
$\Ece{\Cof}\equiv_c \Eceset\equiv_c E^2_=$.
\end{thm}
\begin{pf}
We first show that $\Eceset\leq_c E^2_{=}$. There is a $\Sigma^0_3$ predicate $R(i,x)$ which holds iff $\exists n ( W_x^{[n]}=W_i )$. Let $f(x)$ be a computable function such that $R(i,x)$ iff $i\in W_{f(x)}^{\emptyset''}$. It is then easy to verify that $x~\Eceset~y\Leftrightarrow f(x) ~E^2_= ~f(y)$.

Next we show $E^2_{=}\leq_c \Ece{\Cof}$. There is a single $\Sigma^0_3$ predicate $R$
such that for every $a,x$, we have $a\in W^{\emptyset''}_x\Leftrightarrow R(a,x)$.
Since every $\Sigma^0_3$ set is $1$-reducible to the set $\Cof=\{n: W_n=\dom{\varphi_n}$
is cofinite$\}$, let $g$ be a computable function so that
$a\in W^{\emptyset''}_x\Leftrightarrow W_{g(a,x)}$ is cofinite. Now for each $x$ we produce
the c.e.\ set $W_{f(x)}$ such that for each $a\in\omega$ we have
$W_{f(x)}^{[a]}=\dom{\varphi_{g(a,x)}}$. Hence $f$ is a computable function
witnessing  $E^2_{=}\leq_c \Ece{\Cof}$.

Finally we argue that $\Ece{\Cof}\leq_c \Eceset$. Given a c.e.\ set $A$, and $i,n$,
we let $C(i,n)=[0,i]\cup[i+2,i+M+2]$, where $M$ is the smallest number $\geq n$ such that
$M\not\in A^{[i]}$. Hence the characteristic function of $C(i,n)$ is a string of $i+1$
many 1's, followed by a single 0, and followed by $M+1$ many 1's. Since the least
element not in a c.e.\ set never decreases with time, $C(i,n)$ is uniformly c.e.
Note that if $i\neq i'$ then $C(i,n)\neq C(i',n')$. Now let $D(a,b)=[0,a]\cup[a+2,a+b+1]$.

Now let $\widehat{A}$ be a c.e.\ set having exactly the columns
$\{C(i,n)\mid i,n\in\omega\} \cup\{D(a,b)\mid a,b\in\omega\}$. We verify that $A~E_{\Cof}~B$
iff $\widehat{A}~E_{\text{set}}~\widehat{B}$. Again we write $C(A,i,n), C(B,i,n)$ to distinguish
between the different versions. Suppose that $A~E_{\Cof}~B$. Since $D(a,b)$ appear
as columns in both $\widehat{A}$ and $\widehat{B}$, it suffices to check the $C$ columns.
Fix $C(A,i,n)$. If this is finite then it must equal $D(i,b)$ for some $b$, and so appears as a
column of $\widehat{B}$. If $C(A,i,n)$ is infinite then it is in fact cofinite and so every number
larger than $n$ is eventually enumerated in $A^{[i]}$. Hence $B^{[i]}$ is cofinite and so
$C(B,i,m)$ is cofinite for some $m$. Hence $C(A,i,n)=C(B,i,m)=\omega-\{i+1\}$ appears
as a column of $\widehat{B}$. A symmetric argument works to show that each column
of $\widehat{B}$ appears as a column of $\widehat{A}$.

Now assume that $\widehat{A}~E_{\text{set}}~\widehat{B}$. Fix $i$ such that $A^{[i]}$ is cofinite.
Then $C(A,i,n)=\omega-\{i+1\}$ for some $n$. This is a column of $\widehat{B}$.
Since each $D(a,b)$ is finite $C(A,i,n)=C(B,j,m)$ for some $j$. Clearly $i=j$, which
means that $B^{[i]}$ is cofinite. By a symmetric argument we can conclude
that $A~E_{\Cof}~B$.
\qed
\end{pf}

\begin{thm}
\label{thm:E3Z0}
\Ece{3}$\equiv_c \Zce$.
\end{thm}
\begin{pf}
$\Ece{3}\leq_c\Zce$ was shown in \cite[Prop.\ 3.7]{CHM12}. We now prove $\Zce\leq_c\Ece{3}$. Let $F_s(i,j,n)= \frac{|(W_{i,s}\triangle W_{j,s})\rs{n}|}{n}$. Note that for each $i,j,n$, $F_s(i,j,n)$ changes at most $2n$ times. The triangle inequality holds in this case, that is,  for every $s,x,y,z,n$, we have $F_s(x,z,n)\leq F_s(x,y,n)+F_s(y,z,n)$.

Given $i,j,n,p$ where $i<j<n$ and $p>3$ we describe how to enumerate the finite c.e.\ sets $C_{i,j,n,p}(k)$ for $k\in\omega$. We write $C(k)$ instead of $C_{i,j,n,p}(k)$. For each $k$, $C(k)$ is an initial segment of $\omega$ with at most $n^2(n+1)$ many elements.

If $k\geq n$ we let $C(k)=\emptyset$. We enumerate $C(0),\cdots, C(n-1)$ simultaneously.
Each set starts off being empty, and we assume that $F_0(i,j,n)<2^{-p}$. At each stage there
will be a number $M$ such that $C(i)=[0,M]$, and for every $k < n$, $C(k)=[0,M]$ or $[0,M+1]$.
At stage $s>0$ we act only if $F_s(k_0,k_1,n)$ has changed for some $k_0<k_1<n$. Assume $s$
is such a stage. Suppose $C(i)=[0,M-1]$. We make every $C(k)\supseteq [0,M]$; this is possible as
at the previous stage $C(k)=[0,M-1]$ or $[0,M]$. If $F_s(i,j,n)<2^{-p}$ then do nothing else.
In this case every $C(k)$ is equal to $[0,M]$. Suppose that $F_s(i,j,n)\geq 2^{-p}$. Increase
$C(j)=[0,M+1]$. For each $k\neq i,j$ we need to decide if $C(k)=[0,M]$ or $[0,M+1]$.
%\begin{itemize}\item[(i)] $F_s(i,k,n)<2^{-p-3}$. In this case we say that $k$ wants to follow $i$. Keep $C(k)=[0,M]$.\item[(ii)] $F_s(j,k,n)<2^{-p-3}$.  We say that $k$ wants to follow $j$. Increase $C(k)=[0,M+1]$.\item[(iii)] Neither (i) nor (ii) applies. In this case we say that $k$ has no opinion. Find the least $k'<n$ such that $F_s(k,k',n)<2^{-p-3}$ and $k'$ has an opinion. Keep $C(k)=[0,M]$ if $k'$ wants to follow $i$, and increase $C(k)=[0,M+1]$ if $k'$ wants to follow $j$. If no such $k'$ is found we keep $C(k)=[0,M]$.\end{itemize}

To decide this, consider the graph $G_{i,j,n,p,s}$ with vertices labelled $0,\ldots,n-1$. Vertices $k$ and $k'$ are adjacent iff $F_s(k,k',n)<2^{-(p+k+k'+1)}$, i.e. if $W_k\rs{n}$ and $W_{k'}\rs{n}$ are close and have small Hamming distance. It follows easily from the triangle inequality that $i$ and $j$ must lie in different components. If $k$ is in the same component as $j$ we increase $C(k)=[0,M+1]$ and otherwise keep $C(k)=[0,M]$. This ends the description of the construction.

It is clear that $C_{i,j,n,p}(k)$ is an initial segment of $\omega$ with at most $2n\binom{n}{2}= n^2(n+1)$ many elements. For each $k$, define the set $\widehat{W}_k$ by letting $\widehat{W}_k^{[\langle i,j,p\rangle]}=C_{i,j,j+1,p}(k)\star C_{i,j,j+2,p}(k)\star C_{i,j,j+3,p}(k)\star\cdots$ on column $\langle i,j,p\rangle$, where $i<j$ and $p>3$. Here $C_{i,j,j+1,p}(k)\star C_{i,j,j+2,p}(k)$ denotes the set $X$ where $X(z)=C_{i,j,j+1,p}(k)(z)$ if $z\leq (j+1)^2(j+2)$ and $X(z+(j+1)^2(j+2)+1)=C_{i,j,j+2,p}(k)(z)$. Essentially this concatenates the sets,
with $C_{i,j,j+2,p}(k)$ after the set $C_{i,j,j+1,p}(k)$. The iterated $\star$ operation
is defined the obvious way (and $\star$ is associative). We call the copy of $C_{i,j,n,p}(k)$
in $\widehat{W}_k^{[\langle i,j,p\rangle]}$ the $n^{th}$ block of
$\widehat{W}_k^{[\langle i,j,p\rangle]}$.

We now check that the reduction works. Suppose $W_x~Z_0~W_y$, where $x<y$. Hence we have $\lim\sup _n F(x,y,n)=0$. Fix a column $\langle i,j,p\rangle$. We argue that for almost every $n$, $C_{i,j,n,p}(x)=C_{i,j,n,p}(y)$. There are several cases.
\begin{itemize}
\item[(i)] $\{i,j\}=\{x,y\}$. There exists $n_0>i,j$ such that for every $n\geq n_0$
we have $F(x,y,n)< 2^{-p}$. Hence $C_{i,j,n,p}(x)=C_{i,j,n,p}(y)$ for all large $n$.
\item[(ii)] $|\{i,j\}\cap\{x,y\}|=1$. Assume $i=x$ and $j\neq y$; the other cases
will follow similarly. There exists $n_0>i,j,y$ such that for every $n\geq n_0$
we have $F(x,y,n)< 2^{-(p+x+y+1)}$ and so $x,y$ are adjacent in the graph
$G_{i,j,n,p,s}$ where $s$ is such that $F_s(x,y,n)$ is stable. Since $j$ cannot
be in the same component as $x$, we have  $C_{i,j,n,p}(x)=C_{i,j,n,p}(y)$.
\item[(iii)] $\{i,j\}\cap\{x,y\}=\emptyset$. Similar to (ii). Since $x,y$ are adjacent in the graph $G_{i,j,n,p,s}$ then we must have  $C_{i,j,n,p}(x)=C_{i,j,n,p}(y)$.
\end{itemize}
Hence we conclude that $\widehat{W}_x~E_3~\widehat{W}_y$. Now suppose that
$\widehat{W}_x~E_3~\widehat{W}_y$ for $x<y$. Fix $p>2$ and we have
$\widehat{W}_x^{[\langle x,y,p\rangle]}=^* \widehat{W}_y^{[\langle x,y,p\rangle]}$.
So there is $n_0>y$ such that $C_{x,y,n,p}(x)=C_{x,y,n,p}(y)$ for all $n\geq n_0$.
We clearly cannot have $F(x,y,n)\geq 2^{-p}$ for any $n>n_0$ and so
$\lim\sup_n F(x,y,n)\leq 2^{-p}$. Hence we have $W_x~Z_0~W_y$.
\qed
\end{pf}

\begin{thm}
\label{thm:EsetE3}
$\Eceset \not\leq_c \Ece{3}$.
\end{thm}
\begin{pf}
Suppose there is a computable function witnessing $\Eceset\leq_c \Ece{3}$,
and which maps (the index for) a c.e.\ set $A$ to (the index for) $\widehat{A}$,
so that $A ~E_{\text{set}}~B$ iff $\widehat{A}~E_3~\widehat{B}$. Given (indices for)
c.e.\ sets $A$ and $B$, define
\[         F_s(A,B) =
         \begin{cases}
            \max\{z<x: A(z)\neq B(z)\}, & \text{if $x$ first enters $A\cup B$ at stage $s$,}\\
            \max\{z<s: A(z)\neq B(z)\}, & \text{otherwise.}
         \end{cases}
\]
Here we assume that at each stage $s$ at most one new element is enumerated
in $A\cup B$ at stage $s$. One readily verifies that $F_s(A,B)$ is a total
computable function in the variables involved, with $A=^* B$ iff
$\lim\inf_s F_s(A,B)<\infty$.

We define the c.e.\ sets $A,B$ and $C_0,C_1,\cdots$ by the following. Let $A^{[0]}=\omega$ and for $k>0$ let $A^{[k]}=[0,k-1]$. Let $B^{[k]}=[0,k]$ for every $k$. Finally for each $i$ define
$C_i^{[k]}$ to equal
\[
         \begin{cases}
            [0,j], & \text{if $k=2j+1$,}\\
            \omega, & \text{if $k=2j$ and $\exists^\infty s \left(F_s(\widehat{B}^{[i]},\widehat{C_i}^{[i]})=j\right)$,}\\
            \left[0,\max\{s:F_s(\widehat{B}^{[i]},\widehat{C_i}^{[i]})=j\}\right], & \text{if $k=2j$ and $\forall^\infty s \left(F_s(\widehat{B}^{[i]},\widehat{C_i}^{[i]})\neq j\right)$.}
         \end{cases}
\]
By the recursion theorem we have in advance the indices for $C_0,C_1,\cdots$ so the above definition makes sense. Fix $i$. If $\lim\inf_s F_s(\widehat{B}^{[i]},\widehat{C_i}^{[i]})=\infty$ then every column of $C_i$ is a finite initial segment of $\omega$ and thus we have $C_i~E_{\text{set}}~B$.
By assumption we must have $\widehat{C_i} ~E_3~\widehat{B}$ and thus the two sets agree
(up to finite difference) on every column. In particular $\lim\inf_s F_s(\widehat{B}^{[i]}
\widehat{C_i}^{[i]})<\infty$, a contradiction. Hence we must have $\lim\inf_s F_s(\widehat{B}^{[i]}
\widehat{C_i}^{[i]})=j$ for some $j$. The construction of $C$ ensures that $C_i~E_{\text{set}}~A$
which means that $\widehat{C_i} ~\Ece{3}~\widehat{A}$ and so $\widehat{C_i}^{[i]}
=^* \widehat{A}^{[i]}$. Since $\lim\inf_s F_s(\widehat{B}^{[i]},\widehat{C_i}^{[i]})<\infty$
we in fact have $\widehat{B}^{[i]}=^* \widehat{C_i}^{[i]} =^* \widehat{A}^{[i]}$. Since this must
be true for every $i$ we have $\widehat{B}~E_3~ \widehat{A}$ and so $B~E_{\text{set}}~ A$,
which is clearly false since $B$ has no infinite column.
\qed\end{pf}

The result of Theorem \ref{thm:EsetE3} was something of a surprise.
We were able to see how to give a basic module for a computable
reduction from $\Eceset$ to $\Ece3$, in much the same way that
Proposition 3.9 in \cite{CHM12} serves as a basic module
for Theorem 3.10 there.  In the situation of Theorem
\ref{thm:EsetE3}, we were even able to combine finitely many
of these basic modules, but not all $\omega$-many of them.
The following propositions express this and sharpen
our result.  One the one hand, Propositions \ref{prop:binaryEsetE3}
and \ref{prop:ternaryEsetE3} and the ultimate Theorem \ref{thm:finarity}
show that it really was necessary to build infinitely many sets
to prove Theorem \ref{thm:EsetE3}.  On the other hand,
Theorem \ref{thm:EsetE3} shows that in this case
the proposed basic modules cannot be combined by priority
arguments or any other methods.

\begin{prop}
\label{prop:binaryEsetE3}
There exists a binary reduction from $\Ece{\text{set}}$ to $\Ece3$.
That is, there exist total computable functions $f$ and $g$ such that,
for every $x,y\in\omega$, $x~\Ece{\text{set}}~y$ iff $f(x,y)~\Ece3~g(x,y)$.
\end{prop}
\begin{pf}
We begin with a uniform computable ``chip'' function $h$,
such that, for all $i$ and $j$, $W_i=W_j$ iff
$\exists^\infty s~h(s)=\la i,j\ra$.
Next we show how to define $f$.

First, for every $k\in\omega$, $W_{f(x,y)}$ contains all elements
of every even-numbered column $\omega^{[2k]}$.
To enumerate the elements of $W_{g(x,y)}$ from this column,
we use $h$.  At each stage $s+1$ for which there is some $c$ such
that $h(s)$ is a chip for the sets $W_x^{[k]}$ and $W_y^{[c]}$
(i.e.\ the $k$-th and $c$-th columns of $W_x$ and $W_y$, respectively,
identified effectively by some c.e.\ indices for these sets),
we take it as evidence that these two columns may be equal,
and we find the $c$-th smallest element of
$\overline{W_{g(x,y),s}^{[2k]}}$ and enumerate it into
$W_{g(x,y),s+1}$.

The result is that, if there exists some $c$ such that
$W_x^{[k]}=W_y^{[c]}$, then $W_{g(x,y)}^{[2k]}$
is cofinite, since the $c$-th smallest element of its complement
was added to it infinitely often, each time $W_x^{[k]}$ and
$W_y^{[c]}$ received a chip.  (In the language of these constructions,
the $c$-th marker was moved infinitely many times.)
Therefore $W_{g(x,y)}^{[2k]} =^*\omega = W_{f(x,y)}^{[2k]}$
in this case.  Conversely, if for all $c$ we have
$W_x^{[k]}\neq W_y^{[c]}$, then $W_{g(x,y)}^{[2k]}$
is coinfinite, since for each $c$, the $c$-th marker
was moved only finitely many times, and so
$W_{g(x,y)}^{[2k]} \neq^*\omega = W_{f(x,y)}^{[2k]}$.
Thus $W_{g(x,y)}^{[2k]}=^* W_{f(x,y)}^{[2k]}$
iff there exists $c$ with $W_x^{[k]}=W_y^{[c]}$.

Likewise, $W_{g(x,y)}$ contains all elements of each
odd-numbered column $\omega^{[2k+1]}$, and whenever
$h(s)$ is a chip for $W_y^{[k]}$ and $W_x^{[c]}$,
we adjoin to $W_{f(x,y),s+1}$ the $c$-th smallest
element of the column $\omega^{[2k+1]}$ which is not
already in $W_{f(x,y),s}$.  This process is exactly symmetric
to that given above for the even columns, and the result is that
$W_{f(x,y)}^{[2k]}=^* W_{g(x,y)}^{[2k]}$
iff there exists $c$ with $W_y^{[k]}=W_x^{[c]}$.
So we have established that
$$ x~\Ece{\text{set}}~ y \iff f(x,y) ~\Ece3~ g(x,y)$$
exactly as required.
\qed\end{pf}

\begin{prop}
\label{prop:ternaryEsetE3}
There exists a ternary reduction from $\Ece{\text{set}}$ to $\Ece3$.
That is, there exist total computable functions $f$, $g$, and $h$
 such that, for all $x,y, z\in\omega$:
\begin{align*}
&\text{$x~\Ece{\text{set}}~y$ iff $f(x,y,z)~\Ece3~g(x,y,z)$,}\\
&\text{$y~\Ece{\text{set}}~z$ iff $g(x,y,z)~\Ece3~h(x,y,z)$, and}\\
&\text{$x~\Ece{\text{set}}~z$ iff $f(x,y,z)~\Ece3~h(x,y,z)$.}
\end{align*}
\end{prop}
\begin{pf}
To simplify matters, we lift the notation ``$E_{set}$''
to a partial order $\leqset$, defined on subsets of $\omega$ by:
$$ A\leqset B\iff\text{~every column of $A$ appears as a column in $B$}.$$
So $A~E_{set}~B$ iff $A\leqset B$ and $B\leqset A$.

Again we describe the construction of individual columns
of the sets $W_{f(x,y,z)}$, $W_{g(x,y,z)}$, and $W_{h(x,y,z)}$,
using a uniform chip function for equality on columns.
First, for each pair $\la i,j\ra$, we have a column
designated $L_{ij}^x$, the column where we consider
$x$ on the left for $i$ and $j$.  This means that
we wish to guess, using the chip function, whether
the column $W_x^{[i]}$ occurs as a column in $W_y$,
and also whether it occurs as a column in $W_z$.
We make $W_{f(x,y,z)}$ contain all of this column right away.
For every $c$, we move the $c$-th marker in the column $L_{ij}^x$
in both $W_{g(x,y,z)}$ and $W_{h(x,y,z)}$ whenever
either:
\begin{itemize}
\item
the $c$-th column of $W_y$ receives a chip
saying that it may equal $W_x^{[i]}$; or
\item
the $c$-th column of $W_z$ receives a chip
saying that it may equal $W_x^{[j]}$.
\end{itemize}
Therefore, these columns in $W_{g(x,y,z)}$ and $W_{h(x,y,z)}$
are automatically equal, and they are cofinite
(i.e.\ $=^* W_{f(x,y,z)}$ on this column) iff
either $W_x^{[i]}$ actually does equal some column in $W_y$
or $W_x^{[j]}$ actually does equal some column in $W_z$.

The result, on the columns $L^x_{ij}$ for all $i$ and $j$ collectively,
is the following.
\begin{enumerate}
\item
$W_{g(x,y,z)}$ and $W_{h(x,y,z)}$ are always equal to each other
on these columns.
\item
If $W_x\leqset W_y$, then $W_{f(x,y,z)}$, $W_{g(x,y,z)}$,
and $W_{h(x,y,z)}$ are all cofinite on each of these columns.
\item
If $W_x\leqset W_z$, then again $W_{f(x,y,z)}$, $W_{g(x,y,z)}$,
and $W_{h(x,y,z)}$ are all cofinite on each of these columns.
\item
If there exist $i$ and $j$ such that $W_x^{[i]}$
does not appear as a column in $W_y$ and $W_x^{[j]}$
does not appear as a column in $W_z$, then on that
particular column $L^x_{ij}$, $W_{g(x,y,z)}$
and $W_{h(x,y,z)}$ are coinfinite (and equal),
hence $\neq^* W_{f(x,y,z)}=\omega$.
\end{enumerate}
This explains the name $L^x$:  these columns
collectively ask whether either $W_x\leqset W_y$ or $W_x\leqset W_z$.
We have similar columns $L^y_{ij}$ and $L^z_{ij}$,
for all $i$ and $j$, doing the same operations
with the roles of $x$, $y$, and $z$ permuted.

We also have columns $R^z_{ij}$, for all $i,j\in\omega$,
asking about $W_z$ on the right -- that is, asking whether
either $W_x\leqset W_z$ or $W_y\leqset W_z$.
The procedure here, for a fixed $i$ and $j$, sets both $W_{f(x,y,z)}$
and $W_{g(x,y,z)}$ to contain the entire column
$R^x_{ij}$, and enumerates elements of this column
into $W_{h(x,y,z)}$ using the chip function.
Whenever the column $W_x^{[i]}$ receives a chip
indicating that it may equal $W_z^{[c]}$ for some $c$,
we move the $c$-th marker
in column $R^x_{ij}$ in $W_{h(x,y,z)}$.  Likewise,
whenever the column $W_y^{[j]}$ receives a chip
indicating that it may equal $W_z^{[c]}$ for some $c$,
we move the $c$-th marker in $R^x_{ij}$ in $W_{h(x,y,z)}$.
The result of this construction is that the column
$R^x_{ij}$ in $W_{h(x,y,z)}$ is cofinite
(hence $=^* \omega=W_{f(x,y,z)}=W_{g(x,y,z)}$
on this column) iff at least one of $W_x^{[i]}$ and $W_y^{[j]}$
appears as a column in $W_z$.

Considering the columns $R^z_{ij}$ for all $i$ and $j$ together,
we see that:
\begin{enumerate}
\item
$W_{f(x,y,z)}$ and $W_{g(x,y,z)}$ are always equal to $\omega$
on these columns.
\item
If $W_x\leqset W_z$, then $W_{f(x,y,z)}$, $W_{g(x,y,z)}$,
and $W_{h(x,y,z)}$ are all cofinite on each of these columns.
\item
If $W_y\leqset W_z$, then again $W_{f(x,y,z)}$, $W_{g(x,y,z)}$,
and $W_{h(x,y,z)}$ are all cofinite on each of these columns.
\item
If there exist $i$ and $j$ such that neither $W_x^{[i]}$
nor $W_y^{[j]}$ appears as a column in $W_z$, then on that
particular column $R^z_{ij}$, $W_{h(x,y,z)}$
is coinfinite, hence $\neq^* \omega=W_{f(x,y,z)}=W_{g(x,y,z)}$.
\end{enumerate}

Once again, in addition to the columns $R^z_{ij}$,
we have columns $R^x_{ij}$ and $R^y_{ij}$ for all $i$ and $j$,
on which the same operations take place
with the roles of $x$, $y$, and $z$ permuted.

We claim that the sets $W_{f(x,y,z)}$,
$W_{g(x,y,z)}$, and $W_{h(x,y,z)}$
enumerated by this construction satisfy the proposition.
Consider first the question of whether every column
of $W_x$ appears as a column in $W_z$.
This is addressed by the columns labeled $L^x$
and those labeled $R^z$ (which are exactly
the ones whose construction we described in detail.)
If every column of $W_x$ does indeed appear in $W_z$,
then the outcomes listed there show that all three
of the sets $W_{f(x,y,z)}$, $W_{g(x,y,z)}$, and
$W_{h(x,y,z)}$ are cofinite on every one
of these columns.

On the other hand, suppose some column $W_x^{[i]}$
fails to appear in $W_z$.  Suppose further that
$W_x^{[i]}$ also fails to appear in $W_y$.
Then the column $L^x_{ii}$ has the negative outcome:
on this column, we have
$$W_{f(x,y,z)}\neq^*\omega=W_{g(x,y,z)}=W_{h(x,y,z)}.$$
This shows that $\la f(x,y,z), h(x,y,z)\ra$
(and also $\la f(x,y,z), g(x,y,z)\ra$) fail to lie in $\Ece3$,
which is appropriate, since $\la x,z\ra$ (and $\la x,y\ra$)
were not in $\Ece{\text{set}}$.

The remaining case is that some column $W_x^{[i]}$
fails to appear in $W_z$, but does appear in $W_y$.
In this case, some column $W_y^{[j]}$ (namely, the
copy of $W_x^{[i]}$) fails to appear in $W_z$,
and so the negative outcome on the column $R^z_{ij}$
holds:
$$W_{h(x,y,z)}\neq^*\omega=W_{f(x,y,z)}=W_{g(x,y,z)}.$$
This shows that $\la f(x,y,z), h(x,y,z)\ra$
(and also $\la g(x,y,z), h(x,y,z)\ra$) fail to lie in $\Ece3$,
which is appropriate once again, since $\la x,z\ra$ (and $\la y,z\ra$)
were not in $\Ece{\text{set}}$.

Thus, the situation $W_x\not\leqset W_z$
caused $W_{f(x,y,z)}$ and $W_{h(x,y,z)}$ to differ infinitely
on some column, whereas if $W_x\leqset W_z$, then they
were the same on all of the columns $L^x$ and $R^z$.
Moreover, if they were the same, then $W_{g(x,y,z)}$
was also equal to each of them on these columns.
If they differed infinitely, but $W_x\leqset W_y$, then $W_{g(x,y,z)}$
was equal to $W_{f(x,y,z)}$ on all those columns;
whereas if they differed infinitely and
$W_y\leqset W_z$, then $W_{g(x,y,z)}$ was equal to
$W_{h(x,y,z)}$ on all those columns.

The same holds for each of the other five situations:
for instance, the columns $L^y$ and $R^x$ collectively
give the appropriate outcomes for the question of
whether $W_y\leqset W_x$,
while not causing $W_{h(x,y,z)}$ to differ infinitely
from either $W_{f(x,y,z)}$ or $W_{g(x,y,z)}$
on any of these columns
unless (respectively) $W_z\not\leqset W_x$ or $W_y\leqset W_z$.
Therefore, the requirements of the proposition
are satisfied by this construction.
\qed\end{pf}

\section{Introducing Finitary Reducibility}
\label{sec:introfinitary}

Here we formally begin the study of finitary reducibility,
building on the concepts introduced in Propositions
\ref{prop:binaryEsetE3} and \ref{prop:ternaryEsetE3}.
In Theorem \ref{thm:finarity}, we will sketch the proof
that this construction can be generalized
to any finite arity $n$.  That is, we will show that
$\Ece{\text{set}}$ is $n$-arily reducible to $\Ece3$,
under the following definition.
\begin{defn}
\label{defn:narily}
An equivalence relation $E$ on $\omega$ is \emph{$n$-arily reducible}
to another equivalence relation $F$, written $E\leq_c^n F$,
if there exists a computable total function $f:\omega^n\to\omega^n$
(called an \emph{$n$-ary reduction} from $E$ to $F$)
such that, whenever $f(x_0,\ldots,x_{n-1})=(y_0,\ldots,y_{n-1})$ and $i < j < n$,
we have
$$ x_i~E~x_j\iff y_i~F~y_j.$$
%for all tuples $\xvec=(x_0,\ldots,x_{n-1})$ from $\omega^n$.

If such functions exist uniformly for all $n\in\omega$,
then $E$ is \emph{finitarily reducible} to $F$, written $E\leq_c^{<\omega}F$.
\end{defn}

Often it is simplest to think of the $n$-ary reduction $f$ as a function
$g$ from $\omega^n$ to $\omega^n$, writing $\yvec=g(\xvec)=(f(0,\xvec),\ldots,f(n-1,\xvec))$,
in which case the condition says
$$ (\forall i<n)(\forall j<n)~[x_i~E~x_j\iff y_i~F~y_j].$$
Then a finitary reduction is just a function from $\omega^{<\omega}$
to $\omega^{<\omega}$, mapping $n$-tuples $\xvec$ to $n$-tuples
$\yvec$, with the above property.  The following properties are immediate.
\begin{prop}
\label{prop:basics}
Whenever $E\leq_c^{n+1}F$, we also
have $E\leq_c^n F$.
Finitary reducibility implies all $n$-reducibilities, and
computable reducibility $E\leq_c F$ implies finitary reducibility $E\leq_c^{<\omega}F$.
\end{prop}
\begin{pf}
If $E\leq_c^{n+1}F$ via $h$, then $g(\xvec)=(h(\xvec,0))\res n$
is an $n$-reduction.
If $E\leq_c F$ via $f$, then $(x_0,\ldots,x_{n-1})\mapsto (f(x_0),\ldots,f(x_{n-1}))$
is a finitary reduction.
\qed\end{pf}
Unary reducibility is completely trivial, and binary reducibility
$E\leq^2_c F$ is exactly the same concept as $m$-reducibility on sets
$E\leq_m F$, with $E$ and $F$ viewed as subsets of $\omega$
via a natural pairing function.  For $n>2$, however, we believe
$n$-ary reducibility to be a new concept.
To our knowledge, $\Ece{\text{set}}$ and $\Ece3$ form the
first example of a pair of equivalence relations on $\omega$
proven to be finitarily reducible
but not computably reducible.  A simpler example
appears below in Proposition \ref{prop:maxmin}.

\begin{thm}
\label{thm:finarity}
$\Ece{\text{set}}$ is finitarily reducible to $\Ece3$
(yet $\Ece{\text{set}}\not\leq_c\Ece3$, by Theorem \ref{thm:EsetE3}).
\end{thm}
\begin{pf}
Our proof leans heavily on the details from
Propositions \ref{prop:binaryEsetE3} and \ref{prop:ternaryEsetE3},
and we begin by explaining \ref{prop:ternaryEsetE3}
so as to make clear our generalization.
There the columns $L^x$ can be viewed as a way of asking
whether $X$ has anything else in its equivalence class.
A negative answer, meaning that $W_x\not\leqset W_y$
and $W_x\not\leqset W_z$, clearly implies that
neither $\la x,y\ra$ nor $\la x,z\ra$ lies in $\Ece{\text{set}}$.
%A positive answer, on the other hand, could fail
%to imply the $\leqset$ relations, if $W_y\leqset W_x$,
%for instance.  In Proposition \ref{prop:ternaryEsetE3},
%such other cases were handled by $L^y$ or similar columns.
Here we will give a full argument about the possible
equivalence classes into which $E_{set}$ partitions
the $n$ given c.e.\ sets.

For any fixed $n$, consider each possible partition $P$
of the c.e.\ sets $A_1,\ldots,A_n$ (given by
(arbitrary) indices $m_0,\ldots,m_{n-1}$,
with $A_k=m_{k-1}$) into equivalence classes.
If $P$ is consistent with $E_{set}$ (that is, if every
$E_{set}$-class is contained in some $P$-class), then
for each $i,j$ with $\la A_i,A_j\ra\notin P$,
we have two possible relations:  either
$A_i\not\leqset A_j$ or $A_j\not\leqset A_i$.
We consider every possible conjunction of
one of these possibilities for each such pair $\la i,j\ra$.

We illustrate with an example: suppose $n=5$
and $P$ has classes $\{ A_1,A_2\}$, $\{ A_3,A_4\}$,
and $\{ A_5\}$.  One possible conjunction explaining this situation
is:
\begin{align*}
&A_1\not\leqset A_3~\&~A_1\not\leqset A_4~\&~A_2\not\leqset A_3~\&~
A_2\not\leqset A_4~\&\\
&A_1\not\leqset A_5~\&~A_2\not\leqset A_5~\&~
A_3\not\leqset A_5~\&~A_4\not\leqset A_5.
\end{align*}
Another possibility is:
\begin{align*}
&A_1\not\geqset A_3~\&~A_1\not\leqset A_4~\&~A_2\not\leqset A_3~\&~
A_2\not\geqset A_4~\&\\
&A_1\not\geqset A_5~\&~A_2\not\leqset A_5~\&~
A_3\not\geqset A_5~\&~A_4\not\geqset A_5.
\end{align*}
For this $n$ and $P$ there are $2^8$ such possibilities in all,
since there are $8$ pairs $i<j$ with $\la A_i,A_j\ra\notin P$.
If this $P$ is consistent with $E_{set}$, then at least one
of these $2^8$ possibilities must hold.

Now, for every partition $P$ of $\{ A_1,\ldots,A_n\}$ and for every such
possible conjunction (with $k$ conjuncts, say),
we have an infinite set of columns used in
building the sets $\Ahat_1,\ldots,\Ahat_n$.
These columns correspond to elements of $\omega^k$.
In the second possible conjunction in the example
above, the column for $\la i_1,\ldots,i_k\ra$ corresponds
to the question of whether the following holds.
\begin{align*}
&(\exists c~A_1^{[c]}=A_3^{[i_1]})\text{~or~}
(\exists c~A_1^{[i_2]}=A_4^{[c]})\text{~or~}
(\exists c~A_2^{[i_3]}=A_3^{[c]})\text{~or~}
(\exists c~A_2^{[c]}=A_4^{[i_4]})\text{~or~}\\
&(\exists c~A_1^{[c]}=A_5^{[i_5]})\text{~or~}
(\exists c~A_2^{[i_6]}=A_5^{[c]})\text{~or~}
(\exists c~A_3^{[c]}=A_5^{[i_7]})\text{~or~}
(\exists c~A_4^{[c]}=A_5^{[i_8]}).
\end{align*}
As before, a negative answer implies that $P$ is consistent
with $E_{set}$ on these sets.  Conversely, if $P$ is consistent with $E_{set}$,
then at least one of these $2^8$ disjunctions (in this example)
must fail to hold.

With this framework, the actual construction proceeds
exactly as in Proposition \ref{prop:ternaryEsetE3}.
A uniform chip function guesses whether any of these
eight existential (really $\Sigma^0_3$) statements holds.
If any one does hold, then all sets $\Ahat_i$ are
cofinite in the column for this $P$ and this conjunction
and for $\la i_1,\ldots,i_k\ra$.  If the entire disjunction
(as stated here) is false, then $\Ahat_i=^* \Ahat_j$
on this column iff $\la A_i,A_j\ra\in P$.  So, if $P$
is consistent with $E_{set}$, then we have not caused
$\Ahat_i~E_3~\Ahat_j$ to fail for any $\la i,j\ra$ for which
$A_i~E_{set}~A_j$, but we have caused $\Ahat_i~E_3~\Ahat_j$ to fail
whenever $\la A_i,A_j\ra\notin P$.  (Also, if $P$ is inconsistent
with $E_{set}$, then every disjunction has a positive answer,
so every $\Ahat_i$ is cofinite on each of the relevant columns,
and thus they are all $=^*$ there.)

Of course, one of the finitely many possible equivalence relations
$P$ on $\{ A_1,\ldots,A_n\}$ is actually equal to $E_{set}$ there.
This $P$ shows that, whenever $\la A_i,A_j\ra\notin E_{set}$,
we have $\la\Ahat_i,\Ahat_j\ra\notin E_3$; while the argument
above shows that whenever $A_i~E_{set}~A_j$, neither this $P$
nor any other causes any infinite difference between any of the columns
of $\Ahat_i$ and $\Ahat_j$, leaving $\Ahat_i~E_3~\Ahat_j$.
So we have satisfied the requirements of finitary reducibility,
in a manner entirely independent of $n$ and of the choice of sets
$A_1,\ldots,A_n$.
\qed\end{pf}

A full understanding of this proof reveals that it was essential for
each disjunction to consider every one of the sets $A_1,\ldots,A_n$.
If the disjunction caused $\Ahat_1\neq^*\Ahat_2$ on a particular
column, for example, by making $\Ahat_2$ coinfinite on that
column, then the value of $\Ahat_p$ (for $p>2$)
on that column will be either $\neq^*\Ahat_1$ or $\neq^*\Ahat_2$,
and this decision cannot be made at random.  In fact,
one cannot even just guess from $A_p$
whether or not the relevant column $A_1^{[i]}$ which fails to appear in
$A_2$ appears in $A_p$; in the event that it does not appear,
$\Ahat_p$ may need to be not just coinfinite but actually $=^*\Ahat_2$
on that column.  Since $A_p$ is included in the disjunction
(and in the partition $P$ which generated it),
we have instructions for defining $\Ahat_p$:  either we choose
at the beginning to make it $=\Ahat_1 (=\omega)$ on this column,
or we choose at the beginning to keep it $=\Ahat_2$ there.
The partition $P$ is thus essential as a guide.
For a finite number $n$ of sets, there are only finitely
many $P$ to be considered, but on countably many sets
$A_1,A_2,\ldots$ (such as the collection $W_0,W_1,\ldots$
of all c.e.\ sets), there would be $2^\omega$-many possible
equivalence relations.  Even if we restricted to the $\Pi^0_4$
partitions $P$ (which are the only ones that could equal $\Ece{\text{set}}$),
we would not know, for a given $P$, whether $\Ahat_p$
should be kept equal to $\Ahat_1$ or to $\Ahat_2$,
since a $\Pi^0_4$ relation is too complex to allow effective guessing
about whether it contains $\la 1,p\ra$ or $\la 2,p\ra$.

The concept of $n$-ary reducibility could prove to be a useful measure
of how close two equivalence relations $E$ and $F$ come to being computably
reducible.  The higher the $n$ for which $n$-ary reducibility holds,
the closer they are, with finitary reducibility being the very last step
before actual computable reducibility $E\leq_c F$.
The example of $\Ece{\text{set}}$ and $\Ece3$ is surely quite natural,
and shows that finitary reducibility need not imply computable reducibility.
At the lower levels, we will see in Theorem \ref{thm:3isdifferentfrom4}
that there can also be specific natural differences between
$n$-ary and $(n+1)$-ary reducibility, at least in the case $n=3$.
Another example at the $\Pi^0_2$ level will be given in Proposition \ref{prop:maxmin}.
Right now, though, our first application is to completeness under
these reducibilities.

Working with Ianovski and Nies, we showed in \cite[Thm.\ 3.7 \& Cor.\ 3.8]{IMNNS13}
that no $\Pi^0_{n+2}$ equivalence relation can be complete amongst
all $\Pi^0_{n+2}$ equivalence relations under computable reducibility.
However, we now show that, under finitary reducibility, there is a
complete $\Pi^0_{n+2}$ equivalence relation, for every $n$.
Moreover, the example we give is very naturally defined.
We consider, for each $n$, the equivalence relation $E_{=}^n=\{(i,j)\mid W_i^{\emptyset^{(n)}}=W_j^{\emptyset^{(n)}}\}$. Clearly $E^n_{=}$ is a
$\Pi^0_{n+2}$ equivalence relation. We single out this relation $E^n_=$
because equality amongst c.e.\ sets (and in general, equality amongst
$\Sigma^0_{n+1}$ sets) is indisputably a standard equivalence relation
and, as $n$ varies, permits coding of arbitrary arithmetical information
at the $\Sigma^0_{n+1}$ level.

%Even though $E^n_=$ is not complete amongst $\Pi^0_{n+2}$ equivalence relations
%under computable reducibility, it is complete with respect to the finitary reducibility.
We begin with the case $n=0$.

\begin{thm}
\label{thm:Equalityisuniversal}
The equivalence relation $E^0_=$ (also known as $=^{ce}$)
is complete amongst the $\Pi^0_{2}$ equivalence relations
with respect to the finitary reducibility.
\end{thm}
\begin{pf}
%We first consider $n=0$. The proof for an arbitrary $n>0$ can be obtained by relativization.
Fix a $\Pi^0_2$ equivalence relation $R$. We must produce a computable function $f(k,i,\vec{x})$ such that $f(k,-,-)$ gives the $k$-ary reduction from $R$ to $E^n_=$. Note that the case $k=2$ follows trivially from the fact that $E^0_=$ is $\Pi^0_2$-complete as a set. However the completeness of $E^0_=$ under $\leq^k_c$ for $k>2$ does not follow trivially from this. Nevertheless we will mention the strategy for $k=2$ since it will serve as the basic module.

$k=2$: The strategy for $k=2$ is simple. We monitor the stages at which the pair $(m_0,m_1)$
gets a new chip in $R$. Each time we get a new chip we make $W_{f(2,0,m_0,m_1)}=[0,s]$ and $W_{f(2,1,m_0,m_1)}=[0,s+1]$ where $s$ is a fresh number. Clearly $m_0 R m_1$ iff $W_{f(2,0,m_0,m_1)}=W_{f(2,1,m_0,m_1)}=\omega$. This will serve as the basic module
for the pair $(m_0,m_1)$.

$k=3$: We fix the triple $m_0,m_1,m_2$. For ease of notation we rename these as $0,1,2$ instead. We must build, for $i<3$, the c.e.\ set $A_i=W_{f(3,i,m_0,m_1,m_2)}$. Each $A_i$ will have $\binom{3}{2}=3$ columns, which we denote as $A_i^{a,b}$ for $0\leq a<b< 3$. That is, $A_i^{[0]}=A_i^{0,1},A_i^{[1]}=A_i^{1,2},A_i^{[2]}=A_i^{0,3}$ and $A_i^{[j]}=\emptyset$ for $j>2$.
We assume that at each stage, at most one pair $(i,i')$ gets a new chip.

Each time we get a $(0,1)$-chip we must play the $(0,1)$-game, i.e. we set $A_0^{0,1}=[0,s]$ and $A_1^{0,1}=[0,s+1]$ for a new large number $s$. Of course $A_2^{0,1}$ must decide what to do on this column; for instance if there are infinitely many $(0,2)$-chips then we must make $A_2^{0,1}= A_0^{0,1}$ and if  there are infinitely many $(1,2)$-chips then we must make $A_2^{0,1}= A_1^{0,1}$. At the next stage where we get an $(i,2)$-chip we make $A_2^{0,1}= A_i^{0,1}$. This can be done by padding the shorter column with numbers to match the longer column, and if $A_0^{0,1}$ is made longer then we need to also make $A_1^{0,1}$ longer to keep $A_0^{0,1}\neq A_1^{0,1}$ at every finite stage.

If there are only finitely many $(0,2)$-chips and finitely many $(1,2)$-chips then $\neg 0 R 2$ and $\neg 1 R 2$ and we do not care if $A_2^{0,1}=A_0^{0,1}$ or $A_2^{0,1}=A_1^{0,1}$. Of course $A_2$ has to be different from both $A_0$ and $A_1$ but this will be true at the appropriate columns, i.e. the strategy will ensure that $A_2^{0,2}\neq A_0^{0,2}$ and $A_2^{1,2}\neq A_1^{1,2}$. At some point when the $(i,2)$-chips run out we will stop changing the columns $A_0^{0,1}$ and $A_1^{0,1}$ due to having to ensure the correctness of $A_2$. Hence the outcome of the $(0,1)$-game will be correctly reflected in the columns $A_0^{0,1}$ and $A_1^{0,1}$.

If on the other hand there are infinitely many $(0,2)$-chips and only finitely many $(1,2)$-chips then we have $0R2$ and $\neg 1 R2$. We would have ensured that $A_2^{0,1}=A_0^{0,1}$ (which is important as we must make $A_2=A_0$). Again we do not care if $A_2^{0,1}$ equals $A_1^{0,1}$.

Lastly if there are infinitely many $(i,2)$-chips for each $i<2$ then the interference
of $A_2$ will force both columns $A_0^{0,1}$ and $A_1^{0,1}$ to be $\omega$.
This is acceptable, because $0 R1$ must hold (unless $R$ is not an equivalence relation)
and so the $(0,1)$-game would be played at infinitely many stages anyway.

$k=4$: Again we fix the elements $0,1,2,3$ and build $A_i^{a,b}$ for $i<4$ and $0\leq a<b< 4$. There are now $\binom{4}{2}=6$ columns in each $A_i$. The strategy we used above would seem to suggest in this case that every time we get a $(i,j)$-chip we play the $(i,j)$-game and match columns $A_i^{a,b}$ and $A_j^{a,b}$ whenever $\{a,b\}\cap \{i,j\}=1$. At $n=4$ it is clear that this will not be enough. For instance we could have the equivalence classes $\{0\},\{1\},\{2,3\}$. It could well be that the final $(0,2)$-chip came after the final $(1,2)$-chip, while the final $(1,3)$-chip came after the final $(0,3)$-chip. Then $A_2^{0,1}$ would end up equal to $A_0^{0,1}$ while $A_3^{0,1}$ would end up equal to $A_1^{0,1}$. Since $A_0^{0,1}\neq A_1^{0,1}$ this makes $A_2\neq A_3$, which is not good.

Thus every time $(i,j)$ gets a chip we have to to match columns $A_i^{a,b}$ and $A_j^{a,b}$ for every pair $a,b$ except the pair $(i,j)$. In the above scenario this new rule would force $A_0^{0,1}$ and $A_1^{0,1}$ to increase when a $(2,3)$-chip is obtained. The only way this can happen infinitely often is when $2R3$, and either ($0R2$ and $1R3$) or ($1R2$ and $0R3$). This cycle means that $0R1$ must also be true, and so the $(0,1)$-game would be played infinitely often anyway.

\emph{Arbitrary $k\geq 2$}: We now fix $k\geq 2$, and fix c.e.\ sets $A_0,\ldots,A_{k-1}$.
We describe how to build $A_i^{a,b}$ for $i<k$ and $0\leq a<b<k$. At every stage every
column $A_i^{a,b}$ is just a finite initial segment of $\omega$. We assume at each stage,
at most one chip is obtained. At the beginning enumerate $0$ into $A^{a,b}_b$ for every
$a<b$. At a particular stage in the construction, if no chip is obtained, do nothing.
Otherwise suppose we have an $(i,j)$-chip. We play the $(i,j)$-game, i.e. set  $A_i^{i,j}=[0,s]$
and $A_j^{i,j}=[0,s+1]$ for a fresh number $s$. For each pair $a,b$ such that $(a,b)\neq (i,j)$
we match the columns $A_i^{a,b}$ and $A_j^{a,b}$. What this means is to do nothing if they
are currently equal, and if they are unequal, say $|A_i^{a,b}|<|A_j^{a,b}|$, we fill up $A_i^{a,b}$
with enough numbers to make it equal $A_j^{a,b}$. Furthermore if $a=i$ then $A_b^{a,b}$
should also be topped up to have one more element than $A_i^{a,b}$. This ends the construction
of the columns $A_i^{a,b}$ and of the sets $A_i$.

We now verify that the construction works. It is easy to check that at every stage of the construction, and for every $a<b$ and $i$, we have $|A^{a,b}_a|+1=|A^{a,b}_b|$ and $|A^{a,b}_i|\leq |A^{a,b}_b|$. Now suppose that $iRj$. Then there are infinitely many $(i,j)$-chips obtained during the construction and each time we play the $(i,j)$-game and match every other column of $A_i$ and $A_j$. Hence $A_i=A_j$. Now suppose that $\neg iRj$. We verify that $A^{i,j}_i\neq A^{i,j}_j$. Suppose they are equal, so that they both have to be $\omega$. Let $t_0$ be the stage where the last $(i,j)$-chip is issued. Hence $A^{i,j}_i=[0,s]$ and $A^{i,j}_j=[0,s+1]$ for some fresh number $s$, and so we have $|A^{i,j}_l|\leq |A^{i,j}_i|$ for every $l\neq j$. Let $t_1>t_0$ be the least stage such that either $A^{i,j}_i$ or $A^{i,j}_j$ is increased.

\begin{claim}
If $A_l^{i,j}$ is increased to equal $A_j^{i,j}$ for some $l\neq j$ at some stage $t>t_0$,
then at $t$ some $(l,c)$-chip or $(c,l)$-chip is obtained with $A_c^{i,j}=A_j^{i,j}$.
\end{claim}
\begin{pf}

At $t$ suppose a $(i_0,j_0)$-chip was issued. At $t$ we have three different kind of actions:
\begin{itemize}
\item[(i)] The $(i_0,j_0)$-game is played, affecting columns $A^{i_0,j_0}_{i_0}$ and $A^{i_0,j_0}_{j_0}$.
\item[(ii)] For each $(a,b)\neq (i_0,j_0)$, the smaller of the two columns $A^{a,b}_{i_0}$ or $A^{a,b}_{j_0}$ is increased to match the other.
\item[(iii)] $A^{i_0,b}_b$ is increased in the case $a=i_0$ and $A^{i_0,b}_{i_0}$ is smaller than $A^{i_0,b}_{j_0}$, or $A^{j_0,b}_b$ is increased in the case $a=j_0$ and $A^{j_0,b}_{j_0}$ is smaller than $A^{j_0,b}_{i_0}$.
\end{itemize}
At $t$ the column $A^{i,j}_l$ is increased due to an action of type (i), (ii) or (iii). (i) cannot be because otherwise we have $i_0=i$ and $j_0=j$, but we have assumed that no more $(i,j)$-chips were obtained. It is not possible for (iii) because otherwise $l=j$. Hence we must have (ii) which holds for some $a=i,b=j$. Furthermore $l\in\{i_0,j_0\}$, and letting $c$ be the other element of the set $\{i_0,j_0\}$ we have the statement  of the claim.
\qed\end{pf}

At $t_1$ we cannot have an increase in $A^{i,j}_j$ without an increase in $A^{i,j}_i$, due to the fact that the two always differ by exactly one element. Hence at $t_1$ we know that $A^{i,j}_i$ is increased. It cannot be increased by more than one element because the $(i,j)$-game can no longer be played and we have already seen that $|A^{i,j}_l|\leq |A^{i,j}_j|$ for every $l$. Hence at $t_1$, $A^{i,j}_i$ (and also $A^{i,j}_j$) is increased by exactly one element. Now apply the claim successively to get a sequence of distinct indices $c_0=i,c_1,c_1,c_2,\cdots,c_N=j$ such for every $x$, at least one $(c_x,c_{x+1})$- or $(c_{x+1},c_x)$-chip is obtained in the interval between $t_0$ and $t_1$. Hence we have a new cycle of chips beginning with $i$ and ending with $j$.

Note that at $t_1$, $A^{i,j}_i$ was increased to match $A^{i,j}_c$. Thus the construction at $t_1$ could not have increased the column $A^{i,j}_l$ for any $l\not\in\{i,j\}$. Hence after the action at $t_1$ we again have the similar situation at $t_0$, that is, we again have  $|A^{i,j}_l|\leq |A^{i,j}_i|$ for every $l\neq j$. If $t_1<t_2<t_3<\cdots$ are exactly the stages where $A^{i,j}_i$ or $A^{i,j}_j$ is again increased, we can repeat the claim and the argument above to show that between two such stages we have a new cycle of chips starting with $i$ and ending with $j$. Since there are only finitely many possible cycles, there is a cycle which appears infinitely often, contradicting the transitivity of $R$.

The construction produces a computable function  $f(k,i,\vec{x})$ giving the $k$-ary reduction from the $\Pi^0_2$ relation $R$ to $E^0_=$. Since the construction is uniform in $k$,
finitary reducibility follows.
\qed
\end{pf}

Next we relativize this proof to an oracle.  This will
give $\Pi^0_{n+2}$ equivalence relations which are complete at that level
under finitary reducibility, and will also yield the striking
Corollary \ref{cor:notdcomputable} below, which shows that finitary reductions
can exist even when full reductions of arbitrary complexity fail to exist.
\begin{cor}
\label{cor:Equalityisuniversal}
For each $X\subseteq\omega$, the equivalence relation $E^X_=$ defined by
$$ i~E^X_=~j\qquad\iff\qquad W_i^X = W_j^X$$
is complete amongst all $\Pi^X_{2}$ equivalence relations
with respect to the finitary reducibility.
\end{cor}
\begin{pf}
Essentially, one simply relativizes the entire proof of Theorem
\ref{thm:Equalityisuniversal} to the oracle $X$.  The important point
to be made is that the reduction $f$ thus built is not just $X$-computable,
but actually computable.  Since every set $W_e^X$ in question is now $X$-c.e.,
the program $e=f(i,k,\xvec)$ is allowed to give instructions saying ``look up this
information in the oracle,'' and thus to use an $X$-computable chip function
for an arbitrary $\Pi^X_2$ relation $R$, without actually needing to use
$X$ to determine the program code $e$.
\qed\end{pf}
By setting $X=\emptyset^{(n)}$, we get $\Pi^0_n$-complete equivalence
relations (under finitary reducibility) right up through the arithmetical hierarchy.
\begin{cor}
\label{cor:Pincomplete}
Each equivalence relation $E^n_=$
is complete amongst the $\Pi^0_{n+2}$ equivalence relations
with respect to the finitary reducibility.
\qed\end{cor}

This highlights the central role $E^n_=$ plays amongst the $\Pi^0_{n+2}$ equivalence relations;
it is complete with respect to the finitary reducibility.  A wide variety of
$\Pi^0_{n+2}$ equivalence relations arise naturally in mathematics
(for instance, isomorphism problems for many common classes
of computable structures), and all of these are finitarily reducible to $E_=^n$.
In particular, every $\Pi^0_4$ equivalence relation considered in this section
is finitarily reducible to $E^2_=$.  Indeed, $\Ece{3}$ is complete amongst
$\Pi^0_4$ equivalence relations with respect to the finitary reducibility,
even though $E^2_=\not\leq_c \Ece{3}$.

\begin{cor}
\label{cor:E3complete}
$\Ece{3}$ is complete amongst the $\Pi^0_4$ equivalence relations with respect to the finitary reducibility.
\end{cor}
\begin{pf}
By Theorem \ref{thm:EsetEcof}, $E_=^2~\leq_c~\Eceset$, and
by Theorem \ref{thm:finarity}, $\Eceset~\leq_c^{<\omega}~\Ece{3}$.
The corollary then follows from Corollary \ref{cor:Pincomplete}
and Proposition \ref{prop:basics}.
\qed\end{pf}

Allowing arbitrary oracles in Corollary \ref{cor:Equalityisuniversal}
gives a separate result.  Recall from Definition \ref{defn:compreducibility}
the notion of $\bfd$-computable reducibility.

\begin{cor}
\label{cor:notdcomputable}
For every Turing degree $\bfd$, there exist equivalence relations $E$ and $F$
on $\omega$ such that $E$ is finitarily reducible to $F$ (via a computable function,
of course), but there is no $\bfd$-computable reduction from $E$ to $F$.
\end{cor}
\begin{pf}
We again recall from \cite{IMNNS13} that there is no $\Pi^0_2$-complete
equivalence relation under $\leq_c$.  The proof there relativizes to any
degree $\bfd$ and any set $D\in\bfd$, to show that no $\Pi^D_2$ equivalence
relation on $\omega$ can be complete among $\Pi^D_2$ equivalence
relations even under $\bfd$-computable reducibility.  (The authors of \cite{IMNNS13}
use this relativization to show that there is no $\Pi^0_3$-complete equivalence relation,
for example, by taking $D=\emptyset'$, but their proof really shows that for every
$\Pi^0_3$ equivalence relation, there is another one which is not even
$\bfz'$-computably reducible to the first one.)

Therefore, there exists some $\Pi^D_2$ equivalence relation $E$
such that $E\not\leq_{\bfd}~E^D_=$.  However, Corollary \ref{cor:Equalityisuniversal}
shows that $E$ does have a finitary reduction $f$ to $E^D_=$
(with $f$ specifically shown to be computable, not just $\bfd$-computable).
\qed\end{pf}

\section{Further Results on Finitary Reducibility}
\label{sec:finitary}

\subsection{$\Pi^0_2$ equivalence relations}

Recall the $\Pi^0_2$ \ERs\ $\Ecemin$ and $\Ecemax$, which were defined by
$$ i~\Ecemin~j\iff \min(W_i)=\min(W_j)~~~~~~~~i~\Ecemax~j\iff \max(W_i)=\max(W_j).$$
(Here the empty set has minimum $+\infty$ and maximum $-\infty$,
by definition, while all infinite sets have the same maximum $+\infty$.)
It was shown in \cite{CHM12} that $\Ecemax$ and $\Ecemin$ are both computably
reducible to $\Ece{=}=E^0_=$, and that $\Ecemax$ and $\Ecemin$
are incomparable under $\leq_c$.  The proof given there that
$\Ecemax\not\leq_c\Ecemin$ seemed significantly simpler than the proof
that $\Ecemin\not\leq_c\Ecemax$, but no quantitative distinction
could be expressed at the time to make this intuition concrete.
Now, however, we can use finitary reducibility
to distinguish the two results rigorously.

\begin{prop}\label{prop:maxmin}$\Ecemax$ is not binarily reducible to $\Ecemin$.
However $\Ecemin$ is finitarily reducible to $\Ecemax$.
\end{prop}
\begin{pf}To show $\Ecemax$ is not binarily reducible to $\Ecemin$,
let $f$ be any computable total function.
% we assume that given any pair of c.e.\ sets $W_i,W_j$ we can produce $W_{f(0,i,j)},W_{f(1,i,j)}$ such that $\max W_i=\max W_j$ iff $\min W_{f(0,i,j)} =\min W_{f(1,i,j)}$. %Let $F_s(x,y)=0$ if $\min W_{x,s}=\min W_{y,s}$,
%and $F_s(x,y)=1$ otherwise. Then clearly $F(x,y)=\lim_s F_s(x,y)$ exists, and equals $0$
%iff $\min (W_x)=\min (W_y)$.
%
We build the c.e.\ sets $W_i,W_j$  and assume by the recursion theorem that the indices $i,j$ are given in advance.  At each stage, $W_{i,s}$ and $W_{j,s}$ will both be
initial segments of $\omega$, with $W_{i,0}=W_{j,0}=\emptyset$.
Whenever $\max (W_{i,s})=\max(W_{j,s})$ and $\min(W_{f(0,i,j),s})=\min(W_{f(1,i,j),s})$,
we add the least available element to $W_{i,s+1}$, making the maxima distinct at stage $s+1$.
Whenever $\max (W_{i,s})\neq\max(W_{j,s})$ and $\min(W_{f(0,i,j),s})\neq\min(W_{f(1,i,j),s})$,
we add the least available element to $W_{j,s+1}$, making the maxima the same again.
Since the values of $\min(W_{f(0,i,j),s})$ and $\min(W_{f(1,i,j),s})$ can only change finitely often,
there is some $s$ with $W_i=W_{i,s}$ and $W_j=W_{j,s}$, and our construction
shows that these are both finite initial segments of $\omega$, equal to each other
iff $\min(W_{(f(0,i,j)})\neq\min (W_{f(1,i,j)})$.  Thus $f$ was not a binary reduction.

%Let $s$ be the total number of mind changes by $F_s(f(0,i,j),f(1,i,j))$. Now
%We define $W_i=[0,s]$, by adding the next element each time we see
%$F_{t+1}(f(0,i,j),f(1,i,j)) \neq F_t(f(0,i,j),f(1,i,j))$. Let $W_j=[0,s]$ if $s$ is odd, and $W_j=[0,s)$ if $s>0$ is even, and finally $W_j=\emptyset$ if $s=0$. Now $s<\infty$. If $s=0$ then $F(f(0,i,j),f(1,i,j))=0$ but $W_i\neq\emptyset=W_j$. If $s>0$ is even then again $F(f(0,i,j),f(1,i,j))=0$ but $W_i=[0,s]$ and $W_j=[0,s)$. Finally if $s>0$ is odd then $F(f(0,i,j),f(1,i,j))=1$ but we have $W_i=W_j=[0,s]$, a contradiction.

To show that $\Ecemin$ is finitarily reducible to $\Ecemax$, we must produce a computable function $f(k,i,\vec{x})$ such that $f(k,-,-)$ gives the $k$-ary reduction from $\Ecemin$ to $\Ecemax$. Fixing $k\geq 2$ and indices $m_0,\cdots,m_k$ we describe how to build $W_{f(k,i,\vec{m})}$ for each $i< k$. We denote $A_i=W_{f(k,i,\vec{m})}$. We begin with $A_i=\emptyset$ for all $i$. Each time at a stage $s$ we find a new element enumerated into some $W_{m_i}[s]$ below its current minimum we set $A_j=[0,t+\min W_{m_j}[s]]$ for every $j<k$, where $t$ is a fresh number.

There are only finitely many $m_i$, so $A_j$ is modified only finitely often. So there exists $t$ such that for every $j<k$, $A_j=[0,t+\min W_{m_j}]$. Hence $\min W_{m_i}=\min W_{m_j}$ iff $\max A_i=\max A_j$.
\qed
\end{pf}
This tells us that $\Ecemin\leq_c\Ecemax$ is a lot closer to being true than $\Ecemax\leq_c\Ecemin$.  Surprisingly, we found that the
$\Pi_2^0$ relation $\Ecemax$ is complete for the ternary reducibility
but not for $4$-ary reducibility.

\begin{thm}
\label{thm:3isdifferentfrom4}
$\Ecemax$ is complete for ternary reducibility $\leq_c^3$
among $\Pi^0_2$ equivalence relations, but not so for
$4$-ary reducibility $\leq_c^4$.
\end{thm}
\begin{pf}
By Theorem \ref{thm:Equalityisuniversal}, we may use the relation
$E^0_=$ of equality of c.e.\ sets (also known
as $=^{ce}$), needing only to show that
$E^0_=~\leq_c^3~\Ecemax$ and that
$E^0_=~\not\leq_c^4~\Ecemax$.  First we address
the former claim, building a computable $3$-reduction $f(n,i,j,k)$
as follows.

For any $i,j,k\in\omega$ and any stage $s$, let
$$ m_{ij,s} = \left\{\begin{array}{cl}
s, & \text{if~}W_{i,s}=W_{j,s};\\
\min (W_{i,s}\triangle W_{j,s}), & \text{else.}
\end{array}\right.$$
Thus $W_i\neq W_j$ iff $\lim_s m_{ij,s} <\infty$.
We define $m_{ik,s}$ and $m_{jk,s}$ similarly
for those pairs of sets,
%The construction uses three
%columns from $\omega$:
%$\omega^{[0]}$ to compare $W_i$ and $W_j$,
%$\omega^{[1]}$ to compare $W_i$ and $W_k$,
%and $\omega^{[2]}$ to compare $W_j$ and $W_k$.
%On $\omega^{[0]}$,
%The construction builds three corresponding sets $\What_i$,
%$\What_j$, and $\What_k$, and w
and set $f(0,i,j,k)$,
$f(1,i,j,k)$ and $f(2,i,j,k)$ to be c.e.\ indices of the
three corresponding sets $\What_i$,
$\What_j$, and $\What_k$ built by the following construction.

At each stage $s$, $\What_{i,s}$,
$\What_{j,s}$, and $\What_{k,s}$ will each
be a distinct finite initial segment of $\omega$.  Each time
the sets $W_i$ and $W_j$ get a chip
(i.e.\ appear to be equal), we lengthen each
of these initial segments to be longer than $\What_k$
(but still distinct from each other),
so that $\What_i=\What_j=\omega$
iff $W_i=W_j$, and otherwise they have distinct maxima.
%Meanwhile, $W_{f(2,i,j,k)}$ tries to keep up with
%$W_{f(0,i,j,k)}$ on this column each time $W_k$ and $W_i$
%get a chip, and tries to keep up with $W_{f(1,i,j,k)}$ on it each time
%$W_k$ and $W_j$ get a chip.
Similar arguments apply for $i$ and $k$,
and also for $j$ and $k$.

Let $\What_{i,0}=\{ 0,1\}$, $\What_{j,0}=\{ 0\}$,
and $\What_{k,0}=\emptyset$.
At each stage $s+1$, set $\mhat_s=\max(
\What_{i,s},\What_{j,s},\What_{k,s})$.
We first act on behalf of
$i$ and $j$, checking whether
$m_{ij,s+1}\neq m_{ij,s}$.  If so, then we make $\What_i=[0,\mhat_s+3]$
and $\What_j=[0,\mhat_s+2]$, so that both are longer than
they were before, and if also either $m_{ik,s+1}\neq m_{ik,s}$
or $m_{jk,s+1}\neq m_{jk,s}$, then we set $\What_{k,s+1}=[0,\mhat_s+1]$.
(Otherwise $\What_k$ stays unchanged at this stage.)

If $m_{ij,s+1}=m_{ij,s}$, then we check whether $m_{ik,s+1}\neq m_{ik,s}$.
If so, then we make $\What_i=[0,\mhat_s+3]$
and $\What_k=[0,\mhat_s+2]$, and if also $m_{jk,s+1}\neq m_{jk,s}$,
then we set $\What_{j,s+1}=[0,\mhat_s+1]$.
(Otherwise $\What_j$ stays unchanged at this stage.)

Lastly, if $m_{ij,s+1}=m_{ij,s}$ and $m_{ik,s+1}=m_{ik,s}$,
then we check whether $m_{jk,s+1}\neq m_{jk,s}$.
If so, then we make $\What_j=[0,\mhat_s+3]$
and $\What_k=[0,\mhat_s+2]$, with $\What_i$
staying unchanged.  This completes the construction.

Notice first that if $W_i=W_j$, then $\What_i$
and $\What_j$ were both lengthened at infinitely
many stages, so that $\max(\What_i)=\max(\What_j)=+\infty$.
The same holds for $W_i$ and $W_k$, and also for $W_j$ and $W_k$,
(even though in those cases some of the lengthening may have
come at stages at which we acted on behalf of $W_i$ and $W_j$).
On the other hand, if $W_i\neq W_j$, then at least one of these
must be distinct from $W_k$ as well.  If $W_i\neq W_k$, then
$\What_i$ was lengthened at only finitely many stages;
likewise for $\What_j$ if $W_j\neq W_k$.  So, if two of these sets
were equal but the third was distinct, then the two equal ones
gave rise to sets with maximum $+\infty$ and the third corresponded
to a finite set.  And if all three sets were distinct, then after some stage $s_0$
none of $\What_i$, $\What_j$, and $\What_k$ was ever lengthened again,
in which case they are the three distinct initial segments built at stage $s_0$,
with three distinct (finite) maxima.  So we have defined a ternary reduction
from $E^0_=$ to $\Ecemax$.

However, no $4$-ary relation exists.  We prove this
by a construction using the Recursion Theorem, supposing
that $f$ were a $4$-ary reduction and using indices
$i$, $j$, $k$, and $l$ which ``know their own values.''
We write $\What_i$ for $W_{f(0,i,j,k,l)}$,
$\What_j$ for $W_{f(1,i,j,k,l)}$, and so on as usual,
having first waited for $f$ to converge on these four inputs.
If it converges on them all at stage $s$, we set
$W_{i,s+1}=\{ 0\}$, $W_{j,s+1}=\{ 0,2\}$,
$W_{k,s+1}=\{ 1\}$, and $W_{l,s+1}=\{ 1,3\}$.

Thereafter, at any stage $s+1$ for which
$W_{i,s}\neq W_{j,s}$ and $\max(\What_{i,s})\neq\max(\What_{j,s})$,
we add the next available even number to $W_{i,s+1}$,
leaving $W_{i,s+1}=W_{j,s+1}=W_{j,s}$.
At any stage $s+1$ for which
$W_{i,s}=W_{j,s}$ and $\max(\What_{i,s})=\max(\What_{j,s})$,
we add the next available even number to $W_{j,s+1}$,
leaving $W_{i,s+1}=W_{i,s}\subsetneq W_{j,s+1}$.
Similarly, at any stage $s+1$ for which
$W_{k,s}\neq W_{l,s}$ and $\max(\What_{k,s})\neq\max(\What_{l,s})$,
we add the next available odd number to $W_{k,s+1}$,
leaving $W_{k,s+1}=W_{l,s+1}=W_{l,s}$.
At any stage $s+1$ for which
$W_{k,s}=W_{l,s}$ and $\max(\What_{k,s})=\max(\What_{l,s})$,
we add the next available odd number to $W_{l,s+1}$,
leaving $W_{k,s+1}=W_{l,s}\subsetneq W_{l,s+1}$.
This is the entire construction.

Now if $f$ is indeed a $4$-ary reduction, then it must keep
adding elements to both $\What_i$ and $\What_j$,
since if either of these sets turns out to be finite,
then the construction would have built $W_i$
and $W_j$ to contradict $f$.  So in particular,
$W_i=W_j=\{ 0,2,4,\ldots\}$, and
$\max(\What_i)=\max(\What_j)=+\infty$.  Similarly,
it must keep adding elements to both $\What_k$ and $\What_l$,
and so $W_k=W_l=\{ 1,3,5,\ldots\}$, and
$\max(\What_k)=\max(\What_l)=+\infty$.
But then $W_i\neq W_k$, yet
$\max(\What_i)=\max(\What_k)=+\infty$.
So in fact $f$ was not a $4$-ary reduction.
\qed\end{pf}

The preceding proof of the lack of any $4$-ary reduction
can be viewed as the simple argument that, since
$\Ecemax$ has exactly one $\Pi^0_2$-complete
equivalence class (and all its other classes are $\Delta^0_2$)
while $E^0_=$ has infiinitely many $\Pi^0_2$-complete classes,
the latter cannot reduce to the former.  It requires four
distinct elements of the equivalence relation to show this,
as evidenced by the first half of the proof.
One naturally conjectures that a $\Pi^0_2$ \ER\ with
exactly two $\Pi^0_2$-complete classes might be
complete under $\leq^4_c$, but not under $\leq_c^5$.
In Subsection \ref{subsec:nocollapse} we will see
that this intuition was not correct.

\begin{cor}
\label{cor:relmax}
Theorem \ref{thm:3isdifferentfrom4} relativizes.  That is, for every set $D$,
the equivalence relation $\EDmax$ defined by
$$ i~\EDmax~j\iff \max (W_i^{D})=\max (W_j^D)$$
is complete for ternary computable reducibility $\leq_c^3$
among $\Pi^D_2$ equivalence relations, but not so for
$4$-ary computable reducibility $\leq_c^4$.
\end{cor}
\begin{pf}
Notice that relativizing the proof of Theorem \ref{thm:3isdifferentfrom4}
entirely would give this same result for $D$-computable ternary and
$4$-ary reducibility.  That would be correct, and it follows that $\EDmax$
is not $\Pi^D_2$-complete for $4$-ary computable reducibility $\leq_c^4$
either, since certain $\Pi^D_2$ relations are not even $D$-computable
$4$-arily reducible to it.  However, the ternary completeness
required is also under \emph{computable} reducibility.
Proving it requires the use of the same trick as in Corollary
\ref{cor:Equalityisuniversal}.  Our ternary reduction accepts an input
$\la i,j,k\ra$ and outputs indices $\hat{i}$, $\hat{j}$, and $\hat{k}$
of oracle Turing programs which enumerate $W_i^D$, $W_j^D$, and
$W_j^D$ using their own oracles (since those oracles all happen to be $D$
as well), and then execute the same strategy as in Theorem \ref{thm:3isdifferentfrom4}
for those three sets.
\qed\end{pf}

The relations $\Ecemax$ and $\EDmax$ are quickly seen to be computably
bireducible with the equivalence relations $\Ececard$ and $\EDcard$ (respectively) defined by:
$$ i~\Ececard~j\iff |W_i|=|W_j|~~~~~~~~i~\EDcard~j\iff |W_i^D|=|W_j^D|.$$
So $\Ececard$ is are also $\Pi^0_2$-complete under ternary reducibility
but not under $4$-ary reducibility (by Proposition \ref{prop:basics}),
and similarly with $\EDcard$ for $\Pi^D_2$-completeness under these reducibilities.
The reason for introducing such a similar relation is that a specific relativized
version of it, $\Ecardprime$, appears very useful in computable model theory.
%$$ i~\Ecardprime~j\iff |W^{\emptyset'}_i|=|W^{\emptyset'}_j|.$$
The discussion above, along with Corollary \ref{cor:relmax}, shows
that $\Ecardprime$ is $\Pi^0_3$-complete under ternary $\bfz'$-computable
reducibility but not under $4$-ary computable reducibility.
We will use this fact in the next subsection.
\begin{prop}
\label{prop:Ecardprime}
The equivalence relation $\Ecardprime$ is $\Pi^0_3$-complete
under ternary computable reducibility, but not under $4$-ary
computable reducibility.
\qed\end{prop}

\subsection{Equivalence Relations from Algebra}
\label{subsec:fields}

Having so far considered only equivalence relations from pure computability
theory, we now turn briefly to computable model theory, which one naturally
expects to be a fertile source of equivalence relations.  For background and details
relevant to this section, we refer the reader to \cite{LMS15,M08,M09}.

\begin{defn}
\label{defn:fieldindex}
Fix a computable presentation $K$ of the algebraic closure $\Qbar$
of the rational numbers.  For each $e$, define the field $K_e$ to be the
subfield of $K$ which one gets by closing the c.e.\ subset $W_e$
of the domain of $K$ under the field operations.
The equivalence relation $\FIso$ is now defined to represent
the isomorphism relation among these fields:
$$ i~\FIso~j\iff K_i\cong K_j.$$
\end{defn}

Since every computable algebraic field has a computable embedding into $K$,
with c.e.\ image, we know that the sequence $\la K_e\ra_{e\in\omega}$
includes representatives of every computable algebraic field, up to computable
isomorphism.  Notice also that each $K_e$ may be considered, up to computable
isomorphism, as a computable field itself, since the domain of $K_e$ (which is c.e.,
uniformly in $e$, and infinite) can be pulled back to $\omega$, uniformly in $e$.
In fact, given an $e$ such that $\phi_e$ computes the atomic diagram of a
computable algebraic field of characteristic $0$, one can uniformly find an $i$
and a $j$ such that $\phi_i$ is a computable isomorphism from $K_j$ onto that field.

Algebraically closed fields are usually seen as a simpler class of structures
than algebraic fields, even when the former are allowed to contain transcendental
elements.  In fact, though, the isomorphism problem for computable algebraically
closed fields of characteristic $0$ is $\Pi^0_3$-complete, and thus quantifiably
more difficult than that for computable algebraic fields.  Theorem \ref{thm:FIso}
below shows that the gap
is not as large as suggested by the raw complexity levels:  while $\FIso$
for algebraic fields is $\Pi^0_2$-complete for finitary reducibility, $\FAC$
for computable algebraically closed fields is not $\Pi^0_3$-complete in this way.
Rather, it exhibits the same properties as the relation $\Ecardprime$ from the
preceding subsection:
it is $\Pi^0_3$-complete under ternary reducibility, but not under $4$-ary reducibility.

To make isomorphism on computable algebraically closed fields into an equivalence
relation on $\omega$ in a natural way, we define the field $L_e$
to have transcendence degree $d_e=|\overline{W_e}|$. 
Notice that one can construct a computable copy of this field $L_e$
uniformly effectively in $e$:  for each $n$, we have a field element
$x_n$ which appears to be transcendental over the preceding elements
$x_0,\ldots,x_{n-1}$, but becomes algebraic over $\Q$ if ever $n$ enters $W_e$.
Conversely, given any computable algebraically closed field $F$ of characteristic $0$,
we can find an $i$ with $F\cong L_i$,
effectively in an index $e$ such that $\phi_e$ decides the atomic diagram of $F$.
This is straightforward, since the property of being
algebraically independent over all previous elements of the field is $\Pi^0_1$.
Thus, $L_i\cong L_j$ iff $\overline{W_i}$ and $\overline{W_j}$
have the same size (possibly infinite).  This should immediately
remind the reader of $\Ecardprime$, and indeed, the real content
of the following theorem is that the equivalence relation $\FAC$ defined by
$$ i~\FAC~j\iff L_i\cong L_j$$
is computably bireducible with $\Ecardprime$, while $\FIso$ is bireducible
with $E^0_=$.

\begin{thm}
\label{thm:FIso}
The equivalence relation $\FIso$ on $\omega$, which is $\Pi^0_2$-complete as a set
(under $1$-reducibility), is complete under finitary
reducibility $\leq_c^{<\omega}$ among all $\Pi^0_2$ equivalence relations.
However, the equivalence relation $\FAC$ on $\omega$, which is $\Pi^0_3$-complete
as a set, is only complete under ternary reducibility $\leq_c^3$ among all $\Pi^0_3$ equivalence relations; it is incomplete under $4$-ary reducibility $\leq_c^4$ there.
\end{thm}
\begin{pf}
The $\leq_c^{<\omega}$-completeness result for $\FIso$ follows (using
Proposition \ref{prop:basics} and Theorem \ref{thm:Equalityisuniversal}) from the
computable reduction $f$ from $E_=^0$ to $\leq_c^{<\omega}$ which we now describe.
In $\Qbar$, we can effectively find a square root $\sqrtpn$ of the $n$-th smallest
prime number in the subring $\Z$.  (The two square roots of this prime lie in the same orbit,
so it does not make sense to call them positive or negative.  We simply take the first one we find
in the domain $\omega$ of this presentation of $\Qbar$.)  The key fact here
is that, for each $n$, $\sqrtpn$ does not lie in the subfield generated by the
set $\set{\sqrt{p_m}}{m\neq n}$.  Thus, adjoining any collection of these square
roots to $\Q$ to form a field will not cause any other square roots $\sqrtpn$
to appear in that field.  Therefore, our computable reduction $f$ simply maps each $e$
to an index $f(e)$ of the c.e.\ set $\set{\sqrtpn\in\Qbar}{n\in W_e}$.

On the other hand, we will show that $\FAC$ and $\Ecardprime$ are computably bireducible.
(Recall that $\Ecardprime$ is the relation which holds of indices of $\Sigma^0_2$
sets which have the same cardinality.) Propositions \ref{prop:Ecardprime}
and \ref{prop:basics} then complete our argument.
The computable reduction $h$ from $\FAC$ to $\Ecardprime$ is easy:
just let $W_{h(e)}^{\emptyset'}$ enumerate the elements of $\overline{W_e}$.

For the computable reduction $g$ from $\Ecardprime$ to $\FAC$,
we define $g(e)$ using a fixed total computable chip function $c(e,n,s)$
with $n\in W_e^{\emptyset'}$ iff only finitely many $s$ have $c(e,n,s)=1$.
Build a computable field extension $F$ of $\Q$, starting with elements $x_{n,0}$
(for every $n$) which do not yet satisfy any algebraic relation over $\Q$.
Go through all pairs $\la n,s\ra$ in turn, and whenever we find that
$c(e,n,s)=1$, we make the current $x_{n,k}$ algebraic over $\Q$ (in some way
consistent with the finite portion of the atomic diagram of $F$ enumerated thus far),
and create a new element $x_{n,k+1}$ of $F$ which does not yet satisfy any
algebraic relation over the existing elements.  As we continue, we fill in all
the atomic facts needed to make $F$ into a computable algebraically closed field;
details may be found in \cite{R60}.  Thus, if $n\in W_e^{\emptyset'}$, then
$x_{n,k_n}$ will stay transcendental forever over the preceding elements
(where $k_n$ is the greatest $k$ for which $x_{n,k}$ ever appears in $F$);
while otherwise all $x_{n,k}$ (for every $k$) will eventually be made algebraic.  Thus
$\set{x_{n,k_n}}{n\in W_e^{\emptyset'}}$ is a transcendence basis for $F$,
and so the transcendence degree of $F$ is the cardinality of $W_e^{\emptyset'}$.
We set $g(e)$ to be an index such that $L_{g(e)}$ is isomorphic
to $F$; this index can be found effectively, as remarked above, and clearly
then $g$ is a computable reduction from $\Ecardprime$ to $\FAC$.
\qed\end{pf}

On the other hand, there do exist natural $\Pi^0_3$ isomorphism problems
which are complete under $\leq_c^{<\omega}$ at that level.  The example
we give here is quick, albeit slightly unnatural, in that the field of the equivalence
relation $\EIso0$ is a proper subset of $\omega$.  (An \emph{equivalence structure}
is just an equivalence relation on the domain $\omega$.)  For details, we refer the reader
to \cite{GK02,LMS15}.

\begin{thm}
\label{thm:eqstructs}
The isomorphism problem $\EIso0$ for the class $\E_0$ of computable equivalence
structures with no infinite classes is $\Pi^0_3$-complete under
$\leq_c^{<\omega}$; indeed, $E_=^1$ is computably reducible to $\EIso0$.
\end{thm}
\begin{pf}
For the computable reduction, given an index $i$ of a set $W_i^{\emptyset'}$,
we build a computable equivalence structure $\S$ with domain $\omega$.
$\S$ begins with infinitely many classes of each odd size.  Whenever we see
an initial segment $\sigma\subseteq\emptyset'_s$ of the stage-$s$ approximation
to $\emptyset'$, and an $n\in\omega$ for which $\Phi^\sigma_{e,s}(e)\converges$,
we add a new equivalence class to $\S$, containing $2n+2$ elements.  As long as
this convergence persists at subsequent stages $t>s$, we keep this class this way.
However, if we ever reach a stage $t>s$ with
$\emptyset'_t\res |\sigma|\neq\emptyset'_s\res |\sigma|$, then we add one more element
to this class, giving it an odd number of elements.  In this case, we start searching again
for a new $\sigma$ for which convergence occurs.  This is the entire construction.

It follows that $\S$ has a class of size $2n+2$ iff $n\in W_i^{\emptyset'}$, and that
$\S$ has infinitely many classes of each odd size.  Hence $E_=^1~\leq_c \EIso{0}$
as required.  (Similar constructions show that $E_=^1 \leq_c E_{\cong}^{\alpha}$ for every
$\alpha\leq\omega$, where this is the isomorphism problem for the class $\E_{\alpha}$
of computable equivalence structures with exactly $\alpha$-many infinite classes.)
\qed\end{pf}

We remark that the completeness results about $\FAC$ can readily be seen
also to hold of computable rational vector spaces, which form an extremely
similar class of structures, and could be conjectured to hold for the class
of all computable models of any other strongly minimal theory for which
the independence relation is $\Pi^0_1$ and the spectrum of computable
models of that theory contains all countable models of the theory.
(In all such classes, the isomorphism relation is determined by the
\emph{dimension}, which is the size of a particular subset of the structure,
usually a maximal independent set.)
On the other hand, it would be natural to investigate
other classes for which the isomorphism problem is $\Pi^0_2$, and to determine
whether their isomorphism problems are also $\Pi^0_2$-complete under
finitary reducibility, as in Theorem \ref{thm:FIso}.

\subsection{Distinguishing Finitary Reducibilities}
\label{subsec:nocollapse}

Theorem \ref{thm:3isdifferentfrom4} implies that
$3$-ary and $4$-ary reducibility are distinct notions,
and it is natural to attempt to extend this result to other finitary
reducibilities.  Above we suggested that one way to do so might be
to create $\Pi^0_2$ \ERs\ in which only finitely many
of the equivalence classes are themselves $\Pi^0_2$-complete
as sets.  (We use the class of $\Pi^0_2$-\ERs\ simply because
it is the one we found useful in the preceding subsection.
The same principle could be applied at the $\Pi^0_p$ or other levels,
for any $p$.)  Theorem \ref{thm:infclasses} below
will prove this attempt to be in vain, but the suspicion
that $n$-ary reducibilities are distinct for distinct $n$
turns out to be well-founded, as we will see in Theorem
\ref{thm:4not5}.

It is not difficult to create a $\Pi^0_2$ \ER\ $E$ on $\omega$
having exactly $c$ distinct $\Pi^0_2$-complete equivalence classes.
Define $m~E~n$ iff:
$$(\exists i<m)[m\equiv n\equiv i~(\text{mod}~c)~\&~
\max (W_{\frac{m-i}{c}}) = \max (W_{\frac{n-i}{c}}) ].$$
This essentially just partitions $\omega$ into $c$ distinct classes
modulo $c$, and then partitions each of those classes further
using the relation $\Ecemax$.  As with $\Ecemax$, we intend here
that $\max(W)=\max(V)$ iff $W$ and $V$ are both infinite or both empty
or else have the same (finite) maximum.  For each $i<c$, the class
of those $m\equiv i(\text{mod}~c)$ with $\frac{m-i}{c}\in\textbf{Inf}$
is $\Pi^0_2$-complete, while every other class is defined by such
an $i$ along with a condition of having either a specific finite maximum
(which is a $\Delta^0_1$ condition) or being empty (which is $\Pi^0_1$).

However, this $E$ is not complete among $\Pi^0_2$ \ERs\ under
$4$-ary reducibility.  To build an $F$ with $F\not\leq_c^4 E$,
one uses infinitely many nonconflicting basic modules, one for
each $e$, showing that no $\phi_e$ is a $4$-ary reduction
from $F$ to $E$.  To do this, assign four specific numbers
$w=4e$, $x=4e+1$, $y=4e+2$ and $x=4e+3$ to this module.
Wait until all four of these computations converge:
$\phi_e(1,w,x,y,z)$, $\phi_e(2,w,x,y,z)$, $\phi_e(3,w,x,y,z)$,
and $\phi_e(4,w,x,y,z)$.
(If any diverges, then $\phi_e$ is not total, and we define each of
the four inputs to be an $F$-class unto itself.)  If the four outputs
are all congruent modulo $c$, then we use the same
process which showed that $\Ecemax$ is not $4$-arily
complete for $\Pi^0_2$ \ERs, since now there is
only one $\Pi^0_2$ complete class to which $\phi_e(w)$
and the rest could belong.  On the other hand,
if, say, $\phi_e(1,w,x,y,z)\not\equiv\phi_e(2,w,x,y,z)~(\text{mod}~c)$,
then these two values lie in distinct $E$-classes,
so we just make $w~F~x$; similarly for the other five possibilities.

Nevertheless, there is a straightforward procedure for building an \ER\
which is $4$-complete but not $5$-complete among $\Pi^0_2$ \ERs,
and it generalizes easily to larger finitary reducibilities as well,
showing them all to be distinct.

\begin{thm}
\label{thm:4not5}
For every $n>1$, there exists a $\Pi^0_2$ \ER\ $E$
which is $\Pi^0_2$-complete under $\leq_c^n$, but not under
$\leq_c^{n+1}$.
\end{thm}
\begin{cor}
\label{cor:distinct}
For every $n\neq n'$ in $\omega$, $n$-ary reducibility and
$n'$-ary reducibility do not coincide.
\qed\end{cor}
\begin{pf}
Start with a computable listing $\{( a_{m,0},\ldots,a_{m,n-1})\}_{m\in\omega}$
of all $n$-tuples in $\omega^n$, without repetitions.  The idea is that $E$ should
use the natural numbers $nm,nm+1,\ldots,nm+n-1$ to copy $=^{ce}$
from the $m$-th tuple.
For $i,j\in\omega$, we define $ i~E~j$ if and only if
$$\exists m[
nm\leq i < (n+1)m~~\&~~nm\leq j < (n+1)m
%i,j\in\{ nm, nm+1,\ldots,nm+n-1\}
~~\&~~a_{m,i-mn}=^{ce}a_{m,j-mn}].$$
The last condition just says that $W_{a_{m,i-mn}}=W_{a_{m,j-mn}}$,
which is $\Pi^0_2$.  Of course, for each $i$,
%we can find the only $m\in\omega$ which
only $m=\lfloor\frac{i}{n}\rfloor$ can possibly
satisfy the existential quantifier, %namely that $m$ with $nm\leq i < (n+1)m$,
so this $E$ really is a $\Pi^0_2$ \ER.
Moreover, it is immediate that $=^{ce}$ has an $n$-reduction $f$ to $E$:
for each $n$-tuple $(x_0,\ldots,x_{n-1})\in\omega^n$, just find the unique
$m$ with $(a_{m,0},\ldots,a_{m,n-1}) = (x_0,\ldots,x_{n-1})$, and set
$f(i,x_0,\ldots,x_{n-1})=mn+i$.  That $f$ is an $n$-reduction follows
directly from the design of $E$.  But every $\Pi^0_2$ \ER\ $F$
has an $n$-reduction to $=^{ce}$, since $=^{ce}$ is complete under
finitary reducibility, and so our $E$ is complete under $\leq_c^n$
among $\Pi^0_2$ \ERs.

To show that $E$ is not complete under $\leq_c^{n+1}$, we show that
$=^{ce}~\not\leq_c^{n+1}~E$.  This is surprisingly easy.  Fix any $e\in\omega$,
and define $x_0,\ldots,x_n$ to be the indices of the following programs,
using the Recursion Theorem.  The programs wait until
$\phi_e(i,x_0,\ldots,x_n)$ has converged for every $i\leq n$,
say with $\xhat_i=\phi_e(i,x_0,\ldots,x_n)$.  If all of
$\xhat_0,\ldots,\xhat_n$ lie in a single interval $[nm,(n+1)m)$
for some $m$, then each program $x_i$ simply enumerates $i$
into its set.  Thus we have $x_i\neq^{ce} x_j$ for $i<j\leq n$,
but some two of $\xhat_0,\ldots,\xhat_n$ must be equal,
by the Pigeonhole Principle, and hence $\phi_e$ was not
an $(n+1)$-reduction.  On the other hand, if there exist
$j<k\leq n$ for which $\xhat_j$ and $\xhat_k$ do not lie in the same
interval $[nm,(n+1)m)$, then no program $x_i$ ever enumerates
anything.  In this case we have $x_j =^{ce}x_k$, since both are indices
of the empty set, yet $\la\xhat_j,\xhat_k\ra\notin E$ by the definition of $E$.
Therefore, no $\phi_e$ can be an $(n+1)$-reduction, and so $=^{ce}\not\leq_c^{n+1}E$.
\qed\end{pf}

This proof of Theorem \ref{thm:4not5} is readily adapted to other levels
of the arithmetic hierarchy.  Recall first the following fact.
\begin{prop}
\label{prop:finitarycomplete}
For every $p\geq 0$, there exists a $\Sigma^0_p$ \ER\ which
is complete under finitary reducibility $\leq_c^{<\omega}$ among
$\Sigma^0_p$ \ERs, and a $\Pi^0_p$ \ER\ which is complete
under $\leq_c^{<\omega}$ among $\Pi^0_p$ \ERs.
\end{prop}
\begin{pf}
For $p=0$, equality on $\omega$ is $\Sigma^0_0$-complete
(equivalently, $\Pi^0_0$-complete).  For $p>0$, it is well known
that there is an \ER\
which is $\Sigma^0_p$-complete under full computable reducibility:  let
$\set{V_e}{e\in\omega}$ be a uniform list of the $\Sigma^0_p$ sets,
and take the closure of $\set{(\la e,i\ra,\la e,j\ra)}{\la i,j\ra\in V_e}$
under reflexivity, symmetry, and transitivity.  A $\Pi^0_1$-complete
\ER\ under computable reducibility was constructed in \cite{IMNNS13},
and the \ER\ $\set{(i,j)}{W_i^{\emptyset^{(p-2)}}=W_j^{\emptyset^{(p-2)}}}$
is $\Pi^0_p$-complete under $\leq_c^{<\omega}$ for each $p>1$.
\qed\end{pf}
\begin{thm}
\label{thm:pnotp+1}
For every $p\geq 0$ and every $n\geq 2$, there exists a $\Sigma^0_p$ \ER\ which
is complete under $n$-ary reducibility $\leq_c^{n}$ among
$\Sigma^0_p$ \ERs, but fails to be complete among them
under $\leq_c^{n+1}$.  Likewise, there exists a $\Pi^0_p$ \ER\ which is complete
under $\leq_c^{n}$ among $\Pi^0_p$ \ERs, but not under $\leq_c^{n+1}$.
\end{thm}
\begin{pf}
The $p=0$ case is trivial: every computable \ER\ with
exactly $n$ equivalence classes clearly satisfies the theorem.
Otherwise, the technique is exactly the same as in the proof of Theorem \ref{thm:4not5}.
For $p>0$, fix the $\Sigma^0_p$ \ER\ $F$ which is complete among $\Sigma^0_p$ \ERs\
under $\leq_c^{<\omega}$, as given in Proposition \ref{prop:finitarycomplete}.
Define $ i~E~j$ if and only if
$$\exists m[
nm\leq i < (n+1)m~~\&~~nm\leq j < (n+1)m
~~\&~~a_{m,i-nm}~F~a_{m,j-mn}],$$
again using an effective enumeration
$\set{(a_{m,0},\ldots,a_{m,n-1})}{m\in\omega}$ of $\omega^n$.
Once again we have an $n$-reduction from $F$ to $E$:
set $f(i,x_0,\ldots,x_{n-1})=nm+i$,
where $(a_{m,0},\ldots,a_{m,n-1}) = (x_0,\ldots,x_{n-1})$.
And for $p>0$, the same strategy
as in Theorem \ref{thm:4not5} succeeds in showing that no
$\phi_e$ can be an $(n+1)$-reduction from $F$ to $E$,
although this must be checked for the different cases.
When $p>0$, for each fixed $\phi_e$, there is a computable
reduction to the $\Sigma^0_p$-complete \ER\ $F$ from the
$\Sigma^0_p$ \ER\ which makes $0,\ldots,n$ all equivalent
if all $\phi_e(x_i)$ converge to values in the same
interval $[nm,n(m+1))$, and leaves them pairwise inequivalent otherwise.

The same argument also works with $\Pi^0_p$ in place
of $\Sigma^0_p$.  Our $F$, defined exactly the same way,
is now a $\Pi^0_p$ \ER, and the $n$-ary reduction from $E$
is also the same.  We claim that again $E\not\leq_c^{n+1} F$.
For $p>1$, our $F$ is equality of the sets
$W_i^{\emptyset^{(n)}}$ and $W_j^{\emptyset^{(n)}}$, and so
the proof in Theorem \ref{thm:4not5} using the Recursion Theorem
still works, each c.e.\ set being also c.e.\ in $\emptyset^{(n)}$.
For $p=1$, let all the numbers $\leq n$ be equivalent unless,
on all of those $(n+1)$ numbers, $\phi_e$ converges to values
in the same interval $[nm,n(m+1))$, in which case they become
pairwise inequivalent.  This $\Pi^0_1$ \ER\ must have a computable
reduction to the $\Pi^0_1$-complete \ER\ $F$, which therefore
cannot have any $(n+1)$-ary reduction to $E$.
\qed\end{pf}

Finally, we adapt Theorem \ref{thm:4not5} to compare
finitary reducibility with full computable reducibility.
Of course, it is already known that equality of $\emptyset^{(n)}$-c.e.\
sets is $\Pi^0_{n+2}$-complete under the former, but not
under the latter.
\begin{thm}
\label{thm:finitary}
For each $p>0$,
there exists a $\Sigma^0_p$ \ER\ $E$ which is complete under finitary
reducibility among $\Sigma^0_p$ \ERs, but not under computable reducibility.
\end{thm}
\begin{pf}
Again, let $F$ be $\Sigma^0_p$-complete under computable reducibility.
This time we use an effective enumeration $\{ (a_{m,0},\ldots,a_{m,n_m})\}_{m\in\omega}$
of $\omega^{<\omega}$, and define the computable function
$g$ by $g(0)=\la 0,0\ra$, and
$$g(x+1)=\left\{\begin{array}{cl} \la m,i+1\ra, & \text{if~}g(x)=\la m,i\ra\text{~with~}i<n_m;\\
\la m+1,0\ra, & \textbf{if~}g(x)=\la m,n_m\ra.
\end{array}\right.$$
We let $x~E~y$ iff there is an $m$ with $g(x)=\la m,j\ra$ and
$g(y)=\la m,k\ra$ and $a_{m,j}~F~a_{m,k}$.  Since $F$ is $\Sigma^0_p$,
so is $E$, and the finitary reduction from $F$ to $E$ is given by
$h(i,x_0,\ldots,x_n)=g^{-1}(\la m,i\ra)$, where $(x_0,\ldots,x_n)=(a_{m,0},\ldots,a_{m,n_m})$.
With $F$ $\Sigma^0_p$-complete under $\leq_c$, this makes $E$
$\Sigma^0_p$-complete under $\leq_c^{<\omega}$.
But for each computable total function $f$ (which you think
might be a full computable reduction from $F$ to $E$),
there would be a computable reduction to $E$
from a particular slice of $F$ (say the $c$-th slice)
on which we wait until $f(\la c,0\ra)$ converges to some
number $\la m,k\ra$, then wait until $f$ has converged
on each of $\la c,1\ra,\ldots,\la c,1+n_m\ra$ as well,
and define these $(2+n_m)$ elements to be in distinct $F$-classes
if $f$ maps each of them to a pair of the form $\la m,j\ra$ for the same $m$,
or else all to be in the same $F$-class if not.  As usual,
this shows that $f$ cannot have been a computable reduction.
\qed\end{pf}

So we have answered the basic question.  However, the proof
did not involve any \ER\ with only finitely many $\Pi^0_2$-complete
equivalence classes, as we had originally guessed it would.
Indeed, $4$-completeness for $\Pi^0_2$ \ERs\ turns out to require
a good deal more than just two $\Pi^0_2$-complete equivalence classes,
as we now explain.

Say that a total computable function $h$
is a $\Pi^0_2$-approximating function for an \ER\ $E$ if
$$  (\forall x\forall y)[ x~E~y~~\iff~~\exists^\infty s~h(x,y,s)=1].$$
(We may assume that $h$ has range $\subseteq\{ 0,1\}$.  Every
$\Pi^0_2$ \ER\ has such a function $h$.)
We say that, under this $h$, a particular $E$-class $[z]_E$
is $\Delta^0_2$ if, for all $x,y\in [z]_E$, we have $\lim_s h(x,y,s)=1$.
Of course, if $x\in [z]_E$ and $y\notin [z]_E$, then $\lim_s h(x,y,s)=0$,
so in this case the class $[z]_E$ really is $\Delta^0_2$, uniformly in
any single element $x$ in the class.  On the other hand, even
if $[z]_E$ is not $\Delta^0_2$ under this $h$, it could still
be a $\Delta^0_2$ set, under some other computable approximation.
For this reason, our next theorem does not preclude the possibility that
cofinitely many $E$-equivalence classes might be $\Delta^0_2$,
but it does say that cofinitely many classes cannot be uniformly
limit-computable.

For an example of these notions, let $E$ be the relation $\Ecemax$,
saying of $i$ and $j$ that $W_i$ and $W_j$ have the same maximum.
More formally, $i~\Ecemax~j$ iff
$$ (\forall x\forall s\exists t\exists y,z\geq x)[(x\in W_{i,s}\implies y\in W_{j,t})~\&~
(x\in W_{j,s}\implies z\in W_{i,t})].$$
We can define $h$ here by letting $h(i,j,s) = 1$ when either
$\max(W_{i,s})=\max(W_{j,s})$ or else $\max(W_{i,s})>\max(W_{i,t})$
and $\max(W_{j,s})>\max(W_{j,t})$ (where $t$ is the greatest number $<s$
with $h(i,j,t)=1$), and taking $h(i,j,s)=0$ otherwise.  Then the $\Ecemax$-class
$\textbf{Inf}$ of those $i$ with $W_i$ infinite is the only class which fails
to be $\Delta^0_2$ under this $h$, and since the set $\textbf{Inf}$ is in fact
$\Pi^0_2$-complete, it cannot be $\Delta^0_2$ under any other $h$ either.
Recall that $\Ecemax$ is complete among $\Pi^0_2$ \ERs\ under $\leq_c^3$,
but not under $\leq_c^4$.  The following theorem generalizes this result.

\begin{thm}
\label{thm:infclasses}
Suppose that $E$ is complete under $\leq_c^4$
among $\Pi^0_2$ \ERs.  Let $h$ be any computable
$\Pi^0_2$-approximating function for $E$.  Then
$E$ must contain infinitely many equivalence classes
which are not $\Delta^0_2$ under this $h$.
\end{thm}
\begin{pf}
Suppose that $z_0,\ldots,z_n$ were numbers such that
$\la z_i,z_j\ra\notin E$ for each $i<j$, and such that
%each class $[z_i]_E$ is properly $\Pi^0_2$ under $h$ (i.e. not $\Delta^0_2$ under $h$), but
every $E$-class except these $(n+1)$ classes $[z_i]_E$
is $\Delta^0_2$ under $h$.
For each $e$, we will build four c.e.\ sets
which show that $\phi_e$ is not a $4$-reduction
from the relation $=^{ce}$ to $E$.  (Recall that
$i =^{ce}j$ iff $W_i=W_j$, and that this $\Pi^0_2$-\ER\
is complete under finitary reducibility, making it a natural
choice to show $4$-incompleteness of $E$.)

Fix any $e$, and choose four fresh indices $a$, $b$,
$c$ and $d$ of c.e.\ sets $A=W_a$, $B=W_b$, $C=W_c$,
and $D=W_d$, which we enumerate according to the
following instructions.  First, we wait until $\phi_e(i,a,b,c,d)$
has converged for each $i<4$.  (By the Recursion
Theorem, these indices may be assumed to know their own values.)
Set $\ahat=\phi_e(0,a,b,c,d)$, $\bhat=\phi_e(1,a,b,c,d)$, etc.
If $\phi_e$ is a $4$-reduction, then $A=B$ iff $\ahat~E~\bhat$,
and $A=C$ iff $\ahat~E~\chat$, and so on.

At an odd stage $2s+1$, we first compare $\ahat$ and $\bhat$, using the
computable $\Pi^0_2$-approximating function $h$ for $E$.
If $h(\ahat,\bhat,s)=1$ and $A_{2s}=B_{2s}$, then we
%set $B_{s+1}=A_s$
add to $A_{2s+1}$
some even number not in $B_{2s}$, so $A_{2s+1}\neq B_{2s+1}$.
On the other hand, if $h(\ahat,\bhat,s)=0$ and $A_{2s}\neq B_{2s}$,
then we make $A_{2s+1}=B_{2s+1}=A_{2s}\cup B_{2s}$.  (The purpose
of these maneuvers is to ensure that $\lim_s h(\ahat,\bhat,s)$
diverges, so that $\ahat$ and $\bhat$ lie in one of the properly
$\Pi^0_2$ $E$-classes.)

Next we do exactly the same procedure with $\chat$ and $\dhat$
in place of $\ahat$ and $\bhat$, and using a new odd number if
needed, instead of a new even number.  This completes stage $2s+1$,
ensuring that $\lim_s h(\chat,\dhat,s)$ also diverges.

At stage $2s+2$, fix the $i\leq n$ such that $h(\ahat,z_i,s')=1$ for the greatest
possible $s'\leq s$, and similarly the $j\leq n$ such that
$h(\chat,z_j,s'')=1$ for the greatest possible $s''\leq s$.
(If there are several such $i$, choose the least; likewise for $j$.
If there is no such $i$ or no such $j$, then we do nothing at this stage.)
If $i=j$, then add a new even number to both $A_{2s+2}$ and $B_{2s+2}$,
thus ensuring that they are both distinct from $C_{2s+2}$ and $D_{2s+2}$
(and keeping $A_{2s+2}=B_{2s+2}$ iff $A_{2s+1}=B_{2s+1}$).
If $i\neq j$, then we add all the even numbers in $A_{2s+1}$
to both $C_{2s+2}$ and $D_{2s+2}$, and add all the odd numbers
in $C_{2s+1}$ to both $A_{2s+2}$ and $B_{2s+2}$.
(This is the only step in which even numbers are enumerated
into $C$ or $D$, or odd numbers into $A$ or $B$.)
This completes stage $2s+2$, and the construction.

We claim first that the odd stages succeeded in their purpose
of making $\ahat$, $\bhat$, $\chat$, and $\dhat$ all belong
to properly $\Pi^0_2$ $E$-classes.  At each stage $2s+1$
such that $h(\ahat,\bhat,s)=1$,
we made $A_{2s+1}$ contain a new even number, which only subsequently
entered $B$ if $A_{2s'}=B_{2s'}$ at some stage $s'>s$.
Therefore, if $\lim_s h(\ahat,\bhat,s)=1$, this even number
would show $A\neq B$, yet $\ahat ~E~\bhat$, so that $\phi_e$
would not be a $4$-reduction.  So there are infinitely many $s$
with $h(\ahat,\bhat,s)=0$, and at all corresponding stages
$2s+1$ we made $A_{2s+1}=B_{2s+1}$, which implies
$A=B$.  If $\phi_e$ is a $4$-reduction, then we must have
$\ahat~E~\bhat$, so there were infinitely (but also coinfinitely)
many $s$ with $h(\ahat,\bhat,s)=1$.  Therefore $\lim_s h(\ahat,\bhat,s)$
diverged, and so the $E$-class of $\ahat$ must be one of the $[z_i]_E$
with $i\leq n$, with $\bhat$ lying in the same class.  We now fix this $i$.
A similar analysis on $\chat$ and $\dhat$ shows that they both lie in one
particular $E$-class $[z_j]_E$ with $j\leq n$, and that $C=D$.

Recall that $z_0,\ldots,z_n$ were chosen as representatives of distinct $E$-classes.
Therefore, there must exist some stage $s_0$ such that, at all stages $s>s_0$,
we had $h(\ahat,z_{k},s)=0=h(\bhat,z_{k},s)$ for every $k\neq i$,
and also $h(\chat,z_{k},s)=0=h(\dhat,z_{k},s)$ for every $k\neq j$.
Moreover, we know that $i=j$ iff $z_i~E~z_j$.  If indeed $i=j$,
then at every even stage $>2s_0$ we were in the $i=j$ situation,
and we added a new even number to $A$ and $B$ at each such stage,
while no even numbers were added to either $C$ or $D$ at any
stage $>2s_0$.  Therefore, if $i=j$, we would have $A\neq C$,
yet $\ahat~E~z_i~E~\chat$, which would show that $\phi_e$
is not a $4$-reduction.  On the other hand, if $i\neq j$, then at every even
stage $>2s_0$ we were in the $i\neq j$ situation, and so all even numbers
ever added to $A$ were subsequently added to both $C$ and $D$,
and all odd numbers in $C$ were subsequently added to both $A$ and $B$.
However, no odd numbers were ever added to $A$ or $B$ except numbers
already in $C$, and no even numbers were ever added to $C$ or $D$ except
numbers already in $A$.  So we must have $A=B=C=D$, yet
$\ahat~E~z_i$ and $\chat~E~z_j$, which lie in distinct $E$-classes.
So once again $\phi_e$ cannot have been a $4$-reduction from
$=^{ce}$ to $E$.
This same argument works for every $e$ (by a separate argument for each;
there is no need to combine them), and so $=^{ce}~\!\!\!\not\leq_c^4~E$.
\qed\end{pf}

It remains open whether an \ER\ $E$ which is $\Pi^0_2$-complete under $\leq_c^4$
might have cofinitely many (or possibly all) of its classes be $\Delta^0_2$
in some nonuniform way.

\section{Questions}
\label{sec:questions}

Computable reducibility has been independently invented several times,
but many of its inventions were inspired by the analogy to Borel reducibility
on $2^\omega$.  Therefore, when a new notion appears in computable
reducibility, it is natural to ask whether one can repay some of this debt
by introducing the analogous notion in the Borel context.  We have not attempted
to do so here, but we encourage researchers in Borel reducibility to consider
this idea.  First, do the obvious analogues of $n$-ary and finitary reducibility
bring anything new to the study of Borel reductions?  And second,
in the context of $2^\omega$, could one not also ask about
$\omega$-reducibility?  A Borel $\omega$-reduction from $E$ to $F$
would take an arbitrary countable subset $\{x_0,x_1,\ldots\}$
of $2^\omega$, indexed by naturals, and would produce
corresponding reals $y_0,y_1,\ldots$ with $x_i~E~x_j$ iff
$y_i~F~y_j$.  Obviously, a Borel reduction from $E$ to $F$
immediately gives a Borel $\omega$-reduction, and when the study
of Borel reducibility is restricted to Borel relations on $2^\omega$,
such $\omega$-reductions always exist.  The interesting situation
would involve $E$ and $F$ which are not Borel and for which $E\not\leq_B F$:
could Borel $\omega$-reductions (or finitary reductions) be of use
in such situations?  And finally, if the Continuum Hypothesis fails,
could the same hold true of $\kappa$ reductions, or $<\kappa$-reductions,
for other $\kappa<2^\omega$?

Meanwhile, back on earth, there are plenty of specific questions to be asked
about computable finitary reducibility.  Computable reductions have become
a basic tool in computable model theory, being used to compare classes
of computable structures under the notion of Turing-computable embeddings
(as in \cite{CCKM04,CK06}, for example).  In situations where no
computable reduction exists, finitary reducibility could aid in investigating
the reasons why:  is there not even any binary reduction?  Or is there
a computable finitary reduction, but no computable reduction overall?
Or possibly the truth lies somewhere in between?  Finitary reducibility
has served to answer such questions in several contexts already, as seen
in Subsection \ref{subsec:fields}, and one hopes for it to be used
to sharpen other results as well.

%\bibliographystyle{plain}
%%\begin{singlespace}
%  \bibliography{FinRed, EqRels}
%%\end{singlespace}

\parbox{4.7in}{
{\sc
\noindent
Department of Mathematics \hfill \\
\hspace*{.1in}  Queens College -- C.U.N.Y. \hfill \\
\hspace*{.2in}  65-30 Kissena Blvd. \hfill \\
\hspace*{.3in}  Flushing, New York  11367 U.S.A. \hfill \\
Ph.D. Programs in Mathematics \& Computer Science \hfill \\
\hspace*{.1in}  C.U.N.Y.\ Graduate Center\hfill \\
\hspace*{.2in}  365 Fifth Avenue \hfill \\
\hspace*{.3in}  New York, New York  10016 U.S.A. \hfill}\\
%\medskip
\hspace*{.045in} {\it E-mail: }
\texttt{Russell.Miller\at {qc.cuny.edu} }\hfill \\
\hspace*{.045in} {\it Webpage: }
\texttt{{qcpages.qc.cuny.edu/$\widetilde{~}$rmiller} }\hfill \\
}

\parbox{4.7in}{
{\sc
\noindent
Nanyang Technological University\hfill \\
\hspace*{.1in}  Department of Mathematics \hfill \\
\hspace*{.2in}  Singapore \hfill}\\
%\medskip
\hspace*{.045in} {\it E-mail: }
\texttt{kmng\at {ntu.edu.sg} }\hfill \\
\hspace*{.045in} {\it Webpage: }
\texttt{{www.ntu.edu.sg/home/kmng/} }\hfill \\}

\end{document}